\documentclass[Phys]{SciPost}

\newcommand{\EE}{\mathbb{E}}
\newcommand{\Data}{\mathcal{D}}
\newcommand{\ten}{\mathsf{T}}

\newcommand{\dif}{\partial}

\usepackage{amsmath}
\usepackage{amssymb}
\usepackage{mathtools}
\usepackage{amsthm}
\usepackage{physics}
\usepackage{bm}
\usepackage{comment}
\usepackage{graphicx}
\DeclareMathOperator*{\Extr}{{\rm Extr}}
\newcommand{\oparams}{\vb{\Theta}}

\newcommand{\argmin}{\mathop{\rm arg~min}\limits}
\newtheorem{claim}{Claim}

\newcommand{\usf}{\mathsf{u}}
\newcommand{\vsf}{\mathsf{v}}
\newcommand{\gsf}{\mathsf{g}}
\newcommand{\lsf}{\mathsf{L}}

\newcommand{\tpr}{{t^\prime}}
\newcommand{\spr}{{s^\prime}}

\newcommand{\bhu}{\hat{\bm{u}}}
\newcommand{\bhv}{\hat{\bm{v}}}
\newcommand*{\rev}{\textcolor{black}}

\theoremstyle{plain}

\theoremstyle{definition}

\theoremstyle{remark}

\hypersetup{
    colorlinks,
    linkcolor={red!50!black},
    citecolor={blue!50!black},
    urlcolor={blue!80!black}
}

\DeclareSymbolFont{usualmathcal}{OMS}{cmsy}{m}{n}
\DeclareSymbolFontAlphabet{\mathcal}{usualmathcal}

\fancypagestyle{SPstyle}{
\fancyhf{}
\rhead{{\bf \color{scipostdeepblue} ~Submission }}

\fancyfoot[C]{\textbf{\thepage}}
}

\begin{document}

\pagestyle{SPstyle}

\begin{center}{\Large \textbf{\color{scipostdeepblue}{
%%%%%%%%%% TODO: Write your article's title here
Asymptotic Dynamics of Alternating Minimization for Bilinear Regression\\
%%%%%%%%%% END TODO: TITLE
}}}\end{center}

\begin{center}\textbf{
%%%%%%%%%% TODO: AUTHORS
% Write the author list here. 
% Use (full) first name (+ middle name initials) + surname format.
% Separate subsequent authors by a comma, omit comma and use "and" for the last author.
% Mark the corresponding author(s) with a superscript symbol in this order
% \star, \dagger, \ddagger, \circ, \S, \P, \parallel, ...
Koki Okajima\textsuperscript{1$\star$},
Takashi Takahashi\textsuperscript{1,2$\dagger$}
%%%%%%%%%% END TODO: AUTHORS
}\end{center}

\begin{center}
%%%%%%%%%% TODO: AFFILIATIONS
% Write all affiliations here.
% Format: institute, city, country
{\bf 1} Graduate School of Science, The University of Tokyo, Tokyo, Japan\\
{\bf 2} Institute for Physics of Intelligence, The University of Tokyo, Tokyo, Japan
%%%%%%%%%% END TODO: AFFILIATIONS
%%%%%%%%%% TODO: EMAIL
% Provide email address of corresponding author(s)
\\[\baselineskip]
$\star$ \href{mailto:email1}{\small darjeeling@g.ecc.u-tokyo.ac.jp}\,,\quad
$\dagger$ \href{mailto:email2}{\small takashi-takahashi@g.ecc.u-tokyo.ac.jp}
%%%%%%%%%% END TODO: EMAIL
\end{center}

\section*{\color{scipostdeepblue}{Abstract}}
{\boldmath\textbf{%
%%%%%%%%%% TODO: ABSTRACT
\rev{This study investigates the dynamics of alternating minimization applied 
to a bilinear regression task with normally distributed covariates, 
under the asymptotic system size limit where the number of parameters and observations diverge at the same rate. }
This is achieved by employing the replica method to a multi-temperature glassy system which unfolds the algorithm's time evolution.
Our results show that the dynamics can be described effectively by a two-dimensional discrete stochastic process, 
where each step depends on all previous time steps, revealing the structure of the memory dependence in the evolution of alternating minimization.
The theoretical framework developed in this work can be applied to the analysis of various iterative algorithms, 
extending beyond the scope of alternating minimization.
%%%%%%%%%% END TODO: ABSTRACT
}}

\vspace{\baselineskip}

%%%%%%%%%% BLOCK: Copyright information
% This block will be filled during the proof stage, and finilized just before publication.
% It exists here only as a placeholder, and should not be modified by authors.
%\noindent\textcolor{white!90!black}{%
%\fbox{\parbox{0.975\linewidth}{%
%\textcolor{white!40!black}{\begin{tabular}{lr}%
%  \begin{minipage}{0.6\textwidth}%
%    {\small Copyright attribution to authors. \newline
%    This work is a submission to SciPost Physics. \newline
%    License information to appear upon publication. \newline
%    Publication information to appear upon publication.}
%  \end{minipage} & \begin{minipage}{0.4\textwidth}
%    {\small Received Date \newline Accepted Date \newline Published Date}%
%  \end{minipage}
%\end{tabular}}
%}}
%}
%%%%%%%%%% BLOCK: Copyright information

%%%%%%%%%% TODO: LINENO
% For convenience during refereeing we turn on line numbers:
%\linenumbers
% You should run LaTeX twice in order for the line numbers to appear.
%%%%%%%%%% END TODO: LINENO

%%%%%%%%%% TODO: TOC 
% Guideline: if your paper is longer that 6 pages, include a TOC
% To remove the TOC, simply cut the following block
\vspace{10pt}
\noindent\rule{\textwidth}{1pt}
\tableofcontents
\noindent\rule{\textwidth}{1pt}
%\vspace{10pt}
%%%%%%%%%% END TODO: TOC
%%%%%%%%% TODO: CONTENTS 
% Write your article contents here, starting from first \section.
% An example structure is given below.
\section{Introduction}
\label{sec:intro}
% TODO: write your article here.
Alternating minimization (AM), or classically known as the nonlinear Gauss-Seidel method \cite{Ortega66}, 
is a widely used algorithm for multivariable optimization, where 
one optimizes the objective function with respect to a subset of variables while keeping the rest fixed, 
and then iteratively repeating the process by altering the subset of variables under optimization. 
Explicitly, given the optimization problem 
\begin{equation}
    \min_{\bm{u},\bm{v}\in \mathbb{R}^N} \mathcal{L}(\bm{u},\bm{v}),
\end{equation}
and an initialization point $\bm{u}^0$, the standard AM procedure is given by the following iterative update rule:
\begin{align}\label{eq:AM}
  \begin{split}
    \bhv^{t} &= \argmin_{\bm{v} \in \mathbb{R}^N} \mathcal{L}(\bhu^{t-1},\bm{v}), \\
    \bhu^{t} &= \argmin_{\bm{u} \in \mathbb{R}^N} \mathcal{L}(\bm{u},\bhv^{t}),
  \end{split}
\end{align}
for time indices $t = 1, 2, \cdots .$ 
The use of this technique traces back to the classic Gauss-Seidel method for solving linear systems of equations, 
extending its applications to more contemporary methods such as the EM algorithm \cite{EM77}, matrix factorization \cite{Jain13, Hastie16, Park17}, 
and phase retrieval \cite{Netrapalli13,Waldspurger18,Zhang20}. 
The algorithmic simplicity of AM has been a contributing factor to its popularity 
even for non-convex problems \cite{Wang08, Hu08,Jain13}, in particular to high-dimensional inference tasks, 
where the objective is to retrieve a high-dimensional target vector of size $N$ from a set of $P$ observations.

Although the application of AM is rather ubiquitous, its convergence to a global or satisfactory solution is not guaranteed in general for non-convex problems. 
Therefore, a theoretical understanding of the behavior of such iterative methods is of high interest, as 
it provides insight into their practical effectiveness. 
For example, in inference tasks, the landscape of non-convex objective functions is \rev{known to depend on the relative size of the dataset $P$ compared to the ambient dimension $N$}
\cite{Ge16, Ge17, Ros19, Chi19, Maillard20a, Fyodorov22}. 
In this context, sample complexity analyses have been conducted to determine how much data size $P$ is required to 
accurately retrieve the target vector using AM; a series of studies on low-rank matrix estimation \cite{Jain13,Zhong15} 
and Mixed Linear Regression \cite{Ghosh20} proved 
that AM can recover the target matrix \rev{under sample complexity $ P = O(N \log N)$
using the spectral initialization algorithm given in \cite{YiCaramanis14}. }
Similar results have been obtained for the case of phase retrieval \cite{Waldspurger18}, 
\rev{with necessary sample complexity of $P = O(N)$ for a truncated spectral initialization, }
while $ P = O(N^2)$ for retrieval from a completely \textit{random} initialization. 
Noteworthy progress was made by \cite{Lee23}, proving convergence to the 
target under $O(N\log N)$ sample complexity in rank-one matrix estimation, even under a random initialization.
However, the full characterization of the typical behavior of AM in the asymptotic regime where $P$ and $N$ diverge at the same rate remains an open problem, 
despite there being extensive research in the context of information-theoretic analysis of inference tasks such as phase retrieval \cite{Candes06, Jean19,Maillard20} and low-rank matrix estimation \cite{Jean16,Lelarge17, Barbier24}. 
Our goal is to extend the analysis of AM to this proportional asymptotic regime in order to obtain a sharp characterization of AM under a random design, 
and to provide insights into the impact of initialization and sample complexity that may otherwise be obscured in upper-bound analysis. 
\begin{comment}
While one can appreciate the rigorousness in these studies, 
one may be more interested in understanding the typical behavior of the algorithm under some model assumptions. 
The objective of our analysis is to obtain a sharp characterization of AM 
under a random design, and provide new insight into the impact of initialization and sample complexity 
 which may otherwise be obscured in upper-bound analysis. 
\end{comment}

\rev{
  Statistical physics has provided powerful tools for analyzing the typical behavior of iterative procedures in general, with dynamical mean-field theory (DMFT) being a prominent one. 
  Here, the key idea is to express the generating functional of a dynamical system using path integrals, 
  in which given an iterative procedure dependent on random variables $\bm{J}$, $\bm{x}_{t + 1}= f(\bm{x}_t| \bm{J})$, 
  the average generating functional takes the form 
  \begin{equation}\label{eq:generating_functional_DMFT}
    \EE_{\bm{J}} \mathcal{Z}(\bm{J})[ \qty{\bm{l}_t} ] = \EE_{\bm{J}} \int \prod_{t = 1}^T \dd \bm{x}_t \ \delta\qty( \bm{x}_{t} - f(\bm{x}_{t-1} | \bm{J}) ) e^{\bm{l}_t^\ten \bm{x}_t}.
  \end{equation}
  The tractability of the expectation over $\bm{J}$ crucially depends on the structure of function $f$. 
  For instance, for matrix-valued $\bm{J}$ with i.i.d. Gaussian entries and $f$ being a function of the form $f(\bm{x} | \bm{J}) = f(\bm{J}\bm{x})$, 
  the expectation of this generating functional can be computed by introducing the Fourier representation of the delta function. 
  This approach has enabled the analysis of gradient-based optimization methods \cite{Sarao19,Sarao20,Mignacco20, Sarao21,Mignacco21, Bordelon22, Dandi24,Gerbelot24} and synchronous dynamics of spin glass models \cite{eissfeller92,Opper93,Erba24}. 
  Note that due to the normalization $\mathcal{Z}(\bm{J})[\{ \bm{l}_t = \bm{0} \}] = 1$, it suffices to compute the \textit{annealed} average over $\bm{J}$ to assess the statistical properties of the estimator 
  rather than requiring the \textit{quenched} average over $\bm{J}$. 
  However, as we shall see, the generating functional for AM iterates is not susceptible to this annealed calculation procedure, as the function $f$ given by the argmins \eqref{eq:AM} exhibit a complex dependency on the disorder; 
  see \eqref{eq:online_v}, \eqref{eq:online_u} as well as \eqref{eq:AM_v} and \eqref{eq:AM_u} for a concrete example. 
  Instead, wee express the generating functional as a coupled chain of disordered statistical physics models, where the ground state of each corresponds to the solution of the optimization problem at each iteration. 
  This introduces a non-trivial normalization factor, necessitating a quenched computation over the disorder.  See Section \ref{section:replica_analysis} for a more detailed explanation of our approach. 
}

\rev{The approach of analyzing the dynamics of iterative algorithms as a sequence of glassy systems has been explored in discrete optimization \cite{Krzakala07} and stochastic processes on glassy landscapes \cite{Franz13}. 
However, its full application to optimization algorithms has been limited. 
Recent work has examined two-stage procedures such as knowledge distillation \cite{Saglietti22}, 
and transfer learning \cite{OkajimaObuchi24, Ingrosso24}, where 
the second stage optimization procedure is conditioned on the solution of the first. 
The work by \cite{takahashi22} analyzed the performance of self-training, where the model undergoes an iterative online learning procedure
of creating pseudo-data based on the current model, and then updating the model based on this
generated data. 
 However, a comprehensive analysis of full-batch iterative algorithms under an arbitrary time setup remains unexplored. 
 Our analysis can be seen as an extension of this ``chain of replicas'' approach \cite{Krzakala07} in the context of studying optimization algorithms with full-batched data, 
 which we believe is extendable beyond AM. 
}

\paragraph{Summary of main results.}
In this work, we analyze the dynamics of AM for a bilinear regression task, 
where the objective is to retrieve two target vectors from a set of products of their linear measurement. 
Specifically, our contributions are summarized as follows:
\begin{itemize}
    \item Utilizing the replica method \cite{MPV87,Parisi_Urbani_Zamponi_2020} to compute the quenched average of the generating functional, we provide a closed-form expression 
    for the dynamics of AM (Section \ref{section:replica_analysis}) \rev{in the asymptotic limit where $P, N \to \infty$ with fixed ratio $\kappa := P/N$, 
    while keeping the number of iterations finite in $N$.}
    %This is achieved by unfolding the algorithm's time evolution 
    %into a coupled chain of disordered statistical physics models, each being quenched by the same disorder but holding its own temperature. 
    %The introduces a general framework 
%to study iterative algorithms involving sequential optimization problems. 
    \item Our result offers a statistical characterization of the
    regressors at each iteration by an explicit, discrete two-dimensional Gaussian process, 
    unveiling the memory effect on the algorithm's dynamics \rev{in the effective mean-field picture} (Section \ref{section:Characterization}).
    \item Under this replica analysis, we prove that AM cannot retrieve the target vector 
    under finite $\kappa$ and finite number of iterations when initialized completely randomly ($m_0 = 0$, Subsection \ref{subsection:random_init}). 
    This suggests that fundamentally, the random initialization case 
     requires either a suboptimal number of observations or number of iterations diverging with $N$, 
      which is consistent with previous results for AM with random initialization \cite{Waldspurger18, Lee23} (albeit for different optimization problems).
    \item Comparisons with extensive numerical experiments demonstrate that
    the dynamics for large system size can be captured by our analysis (Section \ref{section:numerical}).
\end{itemize}

\section{The model}
\label{sec:model}

Consider the bilinear regression problem,
 where the objective is to retrieve two unknown target vectors $\bm{u}^\star, \bm{v}^\star \in \mathbb{R}^N$ 
from a dataset $\Data = \qty{ \bm{A}_\mu \in \mathbb{R}^N, \bm{B}_\mu \in \mathbb{R}^N, y_\mu \in \mathbb{R} }_{\mu = 1}^P$, 
with each observation $y_\mu$ given by the product of linear measurements of $\bm{u}^\star$ and $\bm{v}^\star$: 
\begin{equation}
    y_\mu =  (\bm{A}_\mu^\ten \bm{u}^\star )  (\bm{B}_\mu^\ten \bm{v}^\star ),
\end{equation}
 where $\bm{A}_\mu^\ten$ denotes the transpose of $\bm{A}_\mu$ (not to be confused with time index $T$). 
In order to retrieve the target vectors from $\Data$, we consider the following reconstruction scheme via optimization:
\begin{gather}\label{eq:target_problem}
   \min_{\bm{u}, \bm{v}} \mathcal{L}(\bm{u}, \bm{v} | \mathcal{D}), \quad \mathcal{L}(\bm{u}, \bm{v} | \mathcal{D}) =  \sum_{\mu = 1}^P \ell( \bm{A}_\mu^\ten \bm{u}, \bm{B}_\mu^\ten \bm{v};y_\mu ) + \frac{\lambda}{2}\qty( \norm{\bm{u}}_2^2 + \norm{\bm{v}}_2^2 ),
\end{gather}
where $\lambda > 0$ is a regularization parameter. 
Here, the function $\ell(a,b;y)$ is a twice-differentiable bi-convex loss function with respect to $a$ and $b$, and a convex function with respect to $y$.
%While our analysis applies to arbitrary twice-differentiable bi-convex loss functions $l$, we focus on the square loss $l_y(a,b) = (y-ab)^2/2$. 
%In fact, such an optimization task is encountered in several applications such as blind deconvolution \cite{Ahmed14,Xiaodong19,Charisopoulos20} and rank-one matrix sensing \cite{Recht10}. 
Solving for $\bm{u}, \bm{v}$ using AM is a natural approach, as the subproblem at each iteration is essentially a convex optimization problem. 

For the sake of analysis, we assume that the covariates $\{ \bm{A}_\mu, \bm{B}_\mu \}_{\mu = 1}^P$ and target vectors $\bm{u}^\star, \bm{v}^\star$ are
 drawn from an i.i.d. Gaussian ensemble $u_i^\star, v_i^\star \sim \mathcal{N}(0,1) ,\ A_{\mu i}, B_{\mu i} \sim \mathcal{N}(0, 1/N)$ for $\mu = 1, \cdots, P$ and $i = 1, \cdots, N$. 
To investigate the effect of initialization, we also assume 
that the initialization vector of AM, $\bm{u}^0$, has correlation with target $\bm{u}^\star$ controlled by a parameter $m_0$:
\begin{equation}\label{eq:initialization}
    \bm{u}^0 = m_0 \bm{u}^\star + \sqrt{1 - m_0^2} \bm{u}^{\rm n},
\end{equation}
where $\bm{u}^{\rm n}$ is a vector with entries i.i.d. according to $u_i^{\rm n} \sim \mathcal{N}(0, 1), \ i = 1, \cdots , N$. 
%While this initialization setup is purely conventional, we stress that 
%the analysis in our work is extendable to more realistic spectral initializations \cite{Chi19,Xiaodong19,Charisopoulos20}. 
The average with respect to random variables $\{ \Data, \bm{u}^0, \bm{u}^\star, \bm{v}^\star\}$ is denoted by $\EE_\Data$ for brevity. 
Finally, we focus on the high-dimensional setting where the sample complexity is linear with $N$; i.e. $P / N \to \kappa \ (N, P \to \infty)$, 
for $\kappa  = O(1)$.

The objective of our analysis is to precisely determine how the regressors evolve in relation to targets $\bm{u}^\star, \bm{v}^\star$, 
and how much data and good initialization, characterized by the parameters $\kappa$ and $m_0$ respectively, is necessary to retrieve them. 
 In particular, we are interested in the product cosine similarity $m^t$ between the regressors and the targets at any finite iteration $t$ of AM, 
which is defined by 
\begin{equation}\label{eq:product_cosine_similarity}
  m^t := \lim_{N\to \infty} \frac{1}{N} \EE_{\Data} \qty[ \frac{ ( \hat{\bm{u}}^t )^\ten \bm{u}^\star (\hat{\bm{v}}^t)^\ten \bm{v}^\star }{ \norm{\hat{\bm{u}}^t} \norm{\hat{\bm{v}}^t} }].
\end{equation} 
Tracking the evolution of $m^t$ is of particular interest in our analysis, 
as it characterizes the alignment between the estimated and true parameter vectors over the course of the iterative process.

 It should be noted {\color{black} that our problem setup is different from the} \textit{online} setup, where 
 the algorithm is given a new batch of data at each iteration \cite{Zhang20,Verchand23}; i.e. given an 
  initial vector $\bm{u}^0$, the algorithm proceeds as 
  \begin{align}\label{eq:online_v}
    \bhv^{t} &= \argmin_{\bm{v} \in \mathbb{R}^N} \mathcal{L}(\bhu^{t-1},\bm{v} | \mathcal{D}_{t}), \\
    \label{eq:online_u}
    \bhu^{t} &= \argmin_{\bm{u} \in \mathbb{R}^N} \mathcal{L}(\bm{u},\bhv^{t} | \mathcal{D}_{t+1/2}),
  \end{align}
where $\qty{\mathcal{D}_\tau}_{\tau = 1, 3/2, 2, 5/2, \cdots }$ are a sequence of independent data batches. 
Since the regressors at each iteration are only coupled with the previous one, the dynamical analysis becomes much simpler due to its Markovian nature. 
Under such a setting, a sharp characterization of the asymptotics under a random design can be obtained in a rigorous way \cite{Verchand22}.
Their proof is based on leveraging the Convex-Gaussian minimax theorem (CGMT) \cite{Stojnic13, Thrampoulidis18},
 which is a rigorous tool for analyzing the precise asymptotics of
 high-dimensional convex optimization problems under random Gaussian designs, to the analysis of successive optimization procedures \cite{Verchand23}.
While our analysis is non-rigorous, the objective is to provide a similar analysis under a more realistic full-batch setting, in which 
case all iterations are statistically coupled via common data $\mathcal{D}$. 

\section{Replica analysis for alternating minimization}
\label{section:replica_analysis}

\rev{In this section, we describe the key methodology to characterize the dynamical behavior of the vectors $\bm{v}^t, \bm{u}^t$ given by \eqref{eq:AM}. 
As mentioned in the Introduction, DMFT is not directly applicable to the analysis of AM as the updates are given by non-trivial solutions of a series of optimization problems. }
%One might think that the natural approach to treat such dynamical evolution of the system is to use dynamical mean-field theory (DMFT) \cite{Martin73}. 
%However, in the current setup, it is difficult to apply the conventional DMFT straightforwardly because the update rule is given by optimization problems, 
%and the state at step $t$ is not explicitly given as a function of the state at the previous time step; see Subsection \ref{subsection:Relation} for more details. 
Therefore, we will adopt an approach that analyzes the probability density encoding the time evolutions \eqref{eq:AM} as the ground state by using the replica method.

\subsection{Alternating minimization as a stochastic process}
\label{subsection:AM_as_stochastic_process}
Given a fixed set of data $\mathcal{D}$ and target vectors $\bm{u}^\star, \bm{v}^\star $, {\color{black} let us introduce a sequence of 
Boltzmann factors with distinct inverse temperatures $\{\beta_v^{t}, \beta_u^{t}\}_{t\leq T}$:
\begin{align}\label{eq:stochastic_process}
    \begin{split}
    \phi_{\beta_v^{t}}( \bm{v}^t | \bm{u}^{t-1} ) & = \exp\qty[-\beta_v^{t} \mathcal{L} ( \bm{u}^{t-1}, \bm{v}^{t} | \Data) + \beta_v^t \lambda \|\bm{u}^{t-1}\|_2^2 / 2 ]  , \\
    \phi_{\beta_u^{t}}( \bm{u}^t | \bm{v}^{t} ) & = \exp\qty[-\beta_u^{t} \mathcal{L} ( \bm{u}^{t}, \bm{v}^{t} | \Data ) + \beta_u^t \lambda \|\bm{v}^{t}\|_2^2 / 2 ]  ,
    \end{split}
\end{align}
for $t = 1,2, \cdots ,$ where $\bm{u}^0$ is a random vector distributed according to 
\begin{equation}\label{eq:initialization_markov}
    P(\bm{u}^{0} | \bm{u}^\star) = \mathcal{N} ( m_0 \bm{u}^\star, (1 - m_0^2)\vb{I}_N ).
\end{equation}
The crux of our analysis stems from the fact that, by taking the limit $\beta_v^{1} \to \infty,$$ \beta_u^{1} \to \infty,$ 
$ \beta_v^{2} \to \infty, $
$\cdots, $ 
\textit{successively}, 
the joint canonical ensemble given by the Boltzmann factors in \eqref{eq:stochastic_process} converges to the deterministic dynamics precisely given by AM. 
\footnote{
  Note that we have subtracted the terms $\beta_v^t \|\bm{u}^{t-1}\|_2^2 / 2$ and $ \beta_u^t \|\bm{v}^{t}\|_2^2 / 2$ from $\mathcal{L} ( \bm{u}^{t-1}, \bm{v}^{t} | \Data)$ and 
  $ \mathcal{L} ( \bm{u}^{t}, \bm{v}^{t} | \Data ) $ respectively in \eqref{eq:stochastic_process}, as they have no effect on the minimization problem at each iteration of AM given in \eqref{eq:AM}. 
}
More explicitly, given the data $\mathcal{D}$ and initialization $\bm{u}^0$ one can define the following joint distribution of the regressors $\{\bm{u}^{(t)}, \bm{v}^{(t)}\}_{t \leq T}$ :
\footnote{
  We remark that each Boltzmann factor \eqref{eq:stochastic_process} is not normalized; $\int \dd \bm{v}^t \phi_{\beta_v^t} (\bm{v}^t | \bm{u}^{t-1}) \neq 1$, 
  and $\int \dd \bm{u}^t \phi_{\beta_u^t} (\bm{u}^t | \bm{v}^t) \neq 1$. 
This is in contrast to the analysis of the Franz-Parisi potential \cite{FranzParisi95,FranzParisi97,FranzParisi98} and other studies using the same technique \cite{Krzakala07,Franz13}. 
With our construction, the meaning of the joint distribution given by \eqref{eq:joint_distribution} may become ambiguous for finite inverse temperature. 
However, in the limit \eqref{eq:successive_limit}, the measure is still expected to concentrate on the path of the AM algorithm. 
Moreover, this construction can slightly simplify the replica analysis because we do not need to introduce nested replicas as in \cite{Krzakala07,Franz13}.
}
\begin{gather}\label{eq:joint_distribution}
  P\qty( \{ \bm{u}^{t}, \bm{v}^{t} \}_{t \leq T} | \mathcal{D}, \bm{u}^0) = \frac{1}{\mathcal{Z}(\mathcal{D}, \bm{u}^0)} \prod_{t = 1}^T  \phi_{\beta_u^{t}}( \bm{u}^{t} | \bm{v}^{t} ) \phi_{\beta_v^{t}}( \bm{v}^{t} | \bm{u}^{t-1} ),\\
    \label{eq:partition_function}
    \mathcal{Z}(\mathcal{D}, \bm{u}^0) := \prod_{t = 1}^T \int \dd \bm{u}^{t} \dd \bm{v}^{t} \phi_{\beta_u^{t}}( \bm{u}^{t} | \bm{v}^{t} )  \phi_{\beta_v^{t}}( \bm{v}^{t} | \bm{u}^{t-1} ) , 
\end{gather}
A full characterization of this joint canonical ensemble,} thus, indicates that one also obtains a full characterization of AM as well; 
in fact, the average of quantities involving regressors up to iteration $T$ can be assessed by calculating the data and initialization average of the logarithm of the partition function $\mathcal{Z}$ in the limit where 
\begin{equation}\label{eq:successive_limit}
  \lim_{\qty[ \beta_u, \beta_v ] \to \infty} := \lim_{\beta_u^{T} \to \infty} \lim_{\beta_v^{T} \to \infty} \cdots \lim_{\beta_u^{1} \to \infty} \lim_{\beta_v^{1} \to \infty} .
\end{equation}
The problem at hand has been altered from a typical case analysis of an iterative algorithm fed with random data, 
to one of a series of glassy systems coupled to one another, frozen to zero temperature in a successive manner. 

At first glance, the partition function itself seems to be dominated by the contribution of $\beta_v^1$, i.e. $\mathcal{Z}(\mathcal{D}, \bm{u}^0) \simeq \exp \beta_v^{1} \min_{\bm{v}} \mathcal{L} ( \bm{u}^{0}, \bm{v} | \Data )$, which 
may deem the objective of our analysis unachieveable. However, one can consider the logarithm of the partition function as a 
generating function of the regressors $\{\bm{u}^{t}, \bm{v}^{t}\}$. By adding a small external field $\epsilon f$:
\begin{equation}
  \begin{gathered}
  \mathcal{Z}(\mathcal{D}, \bm{u}^0)[\epsilon f(\{\bm{u}^{s}, \bm{v}^{s}\}_{s \leq T})] := \prod_{t = 1}^T \int \dd \bm{u}^{t} \dd \bm{v}^{t}  \phi_{\beta_u^{t}}( \bm{u}^{t} | \bm{v}^{t} ) \phi_{\beta_v^{t}}( \bm{v}^{t} | \bm{u}^{t-1} ) \ e^{\epsilon f(\{\bm{u}^{s}, \bm{v}^{s}\}_{s \leq T})}, 
  \end{gathered}
\end{equation}
one may calculate the average of $f$ under the dynamics of AM by taking the derivative of $\mathcal{Z}(\mathcal{D}, \bm{u}^0)[\epsilon f]$ 
before taking the successive temperature limit: 
{\color{black}
\begin{align}\label{eq:AM_observable}
  \begin{split}
  \langle f \rangle_{\rm AM  | \mathcal{D}, \bm{u}^0 } &= \lim_{\qty[ \beta_u, \beta_v ] \to \infty} \prod_{t =1}^T \int \dd \bm{u}^t \dd \bm{v}^t  P \qty( \{ \bm{u}^{t}, \bm{v}^{t} \} | \mathcal{D}, \bm{u}^0) \  f(\{\bm{u}^{s}, \bm{v}^{s}\}_{s \leq T}) \\
  &=  \lim_{\qty[ \beta_u, \beta_v ] \to \infty} \eval{ \pdv{}{\epsilon}\log \mathcal{Z}(\mathcal{D}, \bm{u}^0)[\epsilon f(\{\bm{u}^{s}, \bm{v}^{s}\}_{s \leq T})]}_{\epsilon = 0}.
  \end{split}
\end{align}}
So far, we have only considered the value of $f$ conditioned on the data $\mathcal{D}$ and initialization $\bm{u}^0$, which is deterministic at this point. 
Here we are interested in the average case analysis with respect to the data and initialization, which accounts to taking the expectation of the right hand side of \eqref{eq:AM_observable} over random data and initialization given in the previous section. 
Assuming that the derivative and expectation can be exchanged, the average generating function can be 
treated using the replica method: 
\begin{equation}\label{eq:replica_method}
    \EE \qty[ \log \mathcal{Z}(\mathcal{D}, \bm{u}^0)] = \lim_{n \to 0} \frac{1}{n} \log \EE\qty[\mathcal{Z}^n], 
\end{equation}
where $\EE$ stands for the joint average over the data and initialization \eqref{eq:initialization_markov}. \rev{
As addressed in the Introduction, recall that the problem is now an analysis of the quenched average of a partition function rather than an annealed average, 
which is what is often encountered in the computational procedure of DMFT.}

\subsection{Outline of the derivation}

In this subsection, we briefly outline the replica computation, i.e., the computation of the right hand side of \eqref{eq:replica_method}. 
See Appendix \ref{appendix:replica_derivation} for the full details of the derivation. 
Readers who are interested in the final expression of the average generating function and its implications may skip this outline and proceed directly to Subsection \ref{subsection:Average_Generating_Function}.

The basic idea of the replica method is to evaluate $\EE[\mathcal{Z}^n]$ for $n\in\mathbb{N}$, and then formally continue the result as $n\to0$ to evaluate the RHS of \eqref{eq:AM_observable}.
Given the statistical properties of the data and initialization given in Section \ref{sec:model}, the $n(\in \mathbb{N})$-th power  of the partition function can be rewritten as 
\begin{equation}\label{eq:replicated_partition_function}
\begin{gathered}
  \EE \qty[\mathcal{Z}^n] = \int \dd \bm{u}^0 \dd \bm{u}^\star \dd \bm{v}^\star P(\bm{u}^0, \bm{u}^\star, \bm{v}^\star) \prod_{a, t =1 }^{n, T} \qty[ \dd \bm{u}^{t}_a \dd \bm{v}^{t}_a  e^{ -\frac{\lambda}{2} ( \beta_u^t \norm{\bm{u}_a^t}_2^2 + \beta_v^t \norm{\bm{v}_a^t}_2^2 )  } ] \\
  \times \Bigg\{ \EE_{\bm{A}, \bm{B}}  \Bigg[ \prod_{a = 1}^n e^{-\beta_v^1 \ell ( \bm{A}^\ten \bm{u}^0 , \bm{B}^\ten \bm{v}_a^1; y ) - \beta_u^1 \ell( \bm{A}^\ten \bm{u}_a^1 , \bm{B}^\ten \bm{v}_a^1 ; y ) } \prod_{t =2 }^{T} e^{  - \beta_v^t \ell( \bm{A}^\ten \bm{u}_a^{t-1}, \bm{B}^\ten \bm{v}_a^t; y ) - \beta_u^t \ell( \bm{A}^\ten \bm{u}_a^t , \bm{B}^\ten \bm{v}_a^t; y ) } \Bigg] \Bigg\}^P ,
\end{gathered}
\end{equation}
where $P(\bm{u}^0, \bm{u}^\star, \bm{v}^\star)$ is the joint distribution of $(\bm{u}^0, \bm{u}^\star, \bm{v}^\star)$, $\bm{A}, \bm{B} \in \mathbb{R}^N$ are 
Gaussian vectors with independent entries of variance $1/N$, and $y = (\bm{A}^\ten \bm{u}^\star)(\bm{B}^\ten \bm{v}^\star)$. 
A crucial observation is that 
the randomness with respect to $\bm{A}, \bm{B}$ only appears via the following random fields : 
\begin{equation}\label{eq:random_fields}
  \begin{gathered}
      h^0 := \bm{A}^\ten \bm{u}^0, \quad h^\star = \bm{A}^\ten\bm{u}^\star,\quad k^\star = \bm{B}^\ten \bm{v}^\star, \quad h_{a}^t := \bm{A}^\ten \bm{u}_a^t, \quad k_{a}^t := \bm{B}^\ten \bm{v}_a^t, \\
       (t = 1, \cdots, T,\ a = 1,\cdots, n).
  \end{gathered}
\end{equation}
Given a fixed configuration of $(\bm{u}^0, \bm{u}^\star, \bm{v}^\star, \{\bm{u}_a^t, \bm{v}_a^t\}_{a,t})$, the random fields are all centered Gaussians with 
covariances given by 
\begin{equation}\label{eq:order_parameters}
  \begin{gathered}
  \EE [ h^\star h^0 ]= \frac{ (\bm{u}^\star)^\ten \bm{u}^0 }{N} =  m^0, \quad \EE [h^t_{a} h^\star] = \frac{ (\bm{u}^t_a)^\ten \bm{u}^\star }{N} =:  m_{u, a}^t, \quad \EE[ k^t_{a} k^\star] = \frac{ (\bm{v}^t_a)^\ten \bm{v}^\star }{N} =:  m_{v, a}^t,   \\
  \EE [ h^t_{a} h^0 ]= \frac{(\bm{u}^t_a)^\ten \bm{u}^0 }{N} =:  R^t_a, \quad \EE [h^s_{a}  h^t_{b}]  =  \frac{ (\bm{u}^s_a)^\ten \bm{u}^t_b }{N} =: Q_{u, ab}^{st}, \quad \EE [k^s_{a} k^t_{b} ]=  \frac{ (\bm{v}^s_a)^\ten \bm{v}^t_b }{N} =: Q_{v, ab}^{st}, \\
  (1 \leq s \leq t \leq T, \ a,b = 1, \cdots, n), 
  \end{gathered}
\end{equation}
where $ \bm{\Theta} = \qty{ m_{u, a}^t, m_{v, a}^t, R^t_a, Q_{u, ab}^{st}, Q_{v, ab}^{st}}$ are the order parameters of the replicated system at hand. 
The order parameters provide the necessary statistics of the Gaussian random fields in \eqref{eq:random_fields}, which simplifies the expression \eqref{eq:replicated_partition_function} to 
\begin{equation}\label{eq:replicated_partition_function_2}
  \begin{gathered}
    \EE \qty[\mathcal{Z}^n] = \int \dd \bm{\Theta} \mathcal{V}(\bm{\Theta})
    %\int \dd \bm{u}^0 \dd \bm{u}^\star \dd \bm{v}^\star P(\bm{u}^0, \bm{u}^\star, \bm{v}^\star) \prod_{a, t =1 }^{n, T} \qty[ \dd \bm{u}^{t}_a \dd \bm{v}^{t}_a  e^{ -\frac{\lambda}{2} ( \beta_u^t \norm{\bm{u}_a^t}_2^2 + \beta_v^t \norm{\bm{v}_a^t}_2^2 )  } ]  \mathcal{V}( \qty{\bm{u}^t_a, \bm{v}^t_a}, \bm{u}^0, \bm{u}^\star, \bm{v}^\star ) \\
 \Bigg\{ \EE \Bigg[ \prod_{a = 1}^n e^{-\beta_v^1 l_{h^\star k^\star} ( h^0, k_a^1 ) - \beta_u^1 l_{h^\star k^\star} ( h_a^1, k_a^1 ) } \prod_{t =2 }^{T} e^{  - \beta_v^t l_{h^\star k^\star} ( h_a^{t-1}, k_a^{t} ) - \beta_u^t l_{h^\star k^\star} ( h_a^t, k_a^t ) } \Bigg] \Bigg\}^P , 
  \end{gathered}
  \end{equation}
where $\mathcal{V}( \bm{\Theta}) $ is the state density of the replicated system satisfying the constraints given in \eqref{eq:order_parameters}, \rev{commonly referred to as the entropic term, 
and the rest corresponding to the energetic term. The entropic term reads }
\begin{align}\label{eq:state_density}
  \begin{split}
  \mathcal{V}(\bm{\Theta}) &= \int \dd \bm{u}^0 \dd \bm{u}^\star \dd \bm{v}^\star P(\bm{u}^0 | \bm{u}^\star) P(\bm{u}^\star) P(\bm{v}^\star) \prod_{a, t =1 }^{n, T} \qty[ \dd \bm{u}^{t}_a \dd \bm{v}^{t}_a  e^{ -\frac{\lambda}{2} ( \beta_u^t \norm{\bm{u}_a^t}_2^2 + \beta_v^t \norm{\bm{v}_a^t}_2^2 )  } ] \\
  &\times  \prod_{a, b = 1}^n \prod_{s \leq t}^T \delta\qty( N Q_{u, ab}^{st} - (\bm{u}^s_a)^\ten \bm{u}^t_b )  \delta\qty( N Q_{v, ab}^{st} - (\bm{v}^s_a)^\ten \bm{v}^t_b ) \\
  & \times \prod_{a,t = 1}^{n,T} \delta\qty( N R^t_a - (\bm{u}^t_a)^\ten \bm{u}^0 ) \delta\qty( N m_{u, a}^t - (\bm{u}^t_a)^\ten \bm{u}^\star ) \delta\qty( N m_{v, a}^t - (\bm{v}^t_a)^\ten \bm{v}^\star ).
 % & = \int \dd \bm{u}^0 \dd \bm{u}^\star \dd \bm{v}^\star P(\bm{u}^0 | \bm{u}^\star) P(\bm{u}^\star) P(\bm{v}^\star) \prod_{a, t =1 }^{n, T} \qty[ \dd \bm{u}^{t}_a \dd \bm{v}^{t}_a  e^{ -\frac{\lambda}{2} ( \beta_u^t \norm{\bm{u}_a^t}_2^2 + \beta_v^t \norm{\bm{v}_a^t}_2^2 )  } ] \\ 
 % &\times \int \dd \hat{\bm{\Theta}} \exp \sum_{(a,b), (s,t)} \qty( N \hat{Q}_{u,ab}^{st} Q_{u,ab}^{st} + N \hat{Q}_{v,ab}^{st} Q_{v,ab}^{st}  - \hat{Q}_u^{ab} (\bm{u}_a^s)^\ten \bm{u}_b^t - \hat{Q}_v^{ab} (\bm{v}_a^s)^\ten \bm{v}_b^t  ) \\
 %&\times \exp N \sum_{a,t} \qty(N \hat{R}_a^t R_a^t + N \hat{m}_{u,a}^t m_{u,a}^t + N \hat{m}_{v,a}^t m_{v,a}^t - (\bm{u}_a^t)^\ten ( \hat{R}_a^t \bm{u}^0 + \hat{m}_{u,a}^t \bm{u}^\star ) - \hat{m}_{v,a}^t (\bm{v}_a^t)^\ten  \bm{v}^\star ),
  \end{split}
\end{align}
%{\color{red} (To be erased : Note that without RS, we still do not know if the inner product between different replicas at different iterations will be symmetric. As a consequence, restricting \textit{both} $a \leq b$ and $s \leq t$ will already be an ansatz. 
%Restricting only $s \leq t$ (or $a \leq b$), on the other hand, will not.)}
In order to obtain an expression that can be continued analytically to $n \to 0$, 
we introduce \textit{replica symmetry} to the set of variables $\vb{\Theta}$, which furthermore constrains the inner products to the following form : 
\begin{equation}
  \begin{gathered}\label{eq:replica_symmetry}
      Q_{u, ab}^{st } = Q_{u,ab}^{ts} = q_{u}^{st} - (1 - \delta_{ab}) \frac{\chi_u^{st}}{\beta^s_u}, \quad Q_{v, ab}^{st } = Q_{v,ab}^{ts} = q_{v}^{st} - (1 - \delta_{ab}) \frac{\chi_v^{st}}{\beta^s_v}  ,\quad 1 \leq s \leq t \leq T, \\
      m_{u, a}^t = m_{u}^t, \quad m_{v,a}^t = m_v^t, \quad R_{a}^t = R^t, \qquad 1 \leq t \leq T. 
  \end{gathered}
\end{equation}
While the validity of replica symmetry is nontrivial, we conjecture that this is true in convex optimization problems, 
based on the experience that replica symmetric computations have been consistent with the other mathematically rigorous analyses \cite{Stojnic13, Thrampoulidis18, Aubin20}. 
 Here, we expect replica symmetry to hold for all iteration index $t$, since each iteration of AM is essentially a minimization of a convex function, albeit being 
 dependent on the solution of the previous one. While we believe that this time-coupling effect does not play a role in replica symmetry breaking, 
 we leave the stability analysis of the replica symmetric solution to future work. 

Given this simplification of order parameters from $\bm{\Theta}$ to $\bm{\Theta}_{\rm RS} := \qty{ q_u^{st}, q_v^{st}, \chi_u^{st}, \chi_v^{st}, m_u^t, m_v^t, R^t }$, 
and by introducing conjugate order parameters $\hat{\bm{\Theta}}_{\rm RS} :=  \qty{ \hat{q}_u^{st}, \hat{q}_v^{st},\hat{\chi}_u^{st}, \hat{\chi}_v^{st}, \hat{m}_u^t, \hat{m}_v^t, \hat{R}^t }$ 
to decouple the delta functions in the entropic term $\mathcal{V}$, 
the Gaussian integrals in the energetic terms in \eqref{eq:replicated_partition_function_2} 
and the high-dimensional integrals in \eqref{eq:state_density} can be further reduced. 

In generic form, the replicated partition function \eqref{eq:replicated_partition_function_2} can be expressed as 
\begin{align}
  \begin{split}
  \EE \qty[ \mathcal{Z}^n ]&= \int \dd \bm{\Theta}_{\rm RS} \dd \hat{\bm{\Theta}}_{\rm RS} \exp nN \qty[ \mathcal{G}(\bm{\Theta}_{\rm RS}, \hat{\bm{\Theta}}_{\rm RS}) + \mathcal{O}(n) ]\\
  &\stackrel{N \to \infty}{\simeq} \exp nN \qty[ \Extr_{\bm{\Theta}_{\rm RS}, \hat{\bm{\Theta}}_{\rm RS}} \mathcal{G}(\bm{\Theta}_{\rm RS}, \hat{\bm{\Theta}}_{\rm RS}) + \mathcal{O}(n) ],
  \end{split}
\end{align}
where we have used the saddle point approximation for large $N$, and $\Extr_x f(x)$ represents the value of $f(x)$ evaluated at its extremum. 
The specific form of the function $\mathcal{G}(\bm{\Theta}_{\rm RS}, \hat{\bm{\Theta}}_{\rm RS})$ is given in the next subsection.
This yields \eqref{eq:replica_method} as an extremum value problem:
\begin{align}
 \lim_{N \to \infty} \frac{1}{N} \EE\qty[ \log \mathcal{Z}(\Data, \bm{u}^0) ]= \Extr_{\bm{\Theta}_{\rm RS}, \hat{\bm{\Theta}}_{\rm RS}} \mathcal{G}(\bm{\Theta}_{\rm RS}, \hat{\bm{\Theta}}_{\rm RS}).
\end{align}

\subsection{Average generating function and saddle point equation}
\label{subsection:Average_Generating_Function}

To provide further detail on the form of the function $\mathcal{G}(\bm{\Theta}_{\rm RS}, \hat{\bm{\Theta}}_{\rm RS})$, 
for convenience we define the following set of order parameters for each time iteration $t$:
\begin{align}
  \bm{\theta}_u^t &:= \qty{ m_u^t,  R^t, \qty{q_u^{st}, \chi_u^{st}}_{s \leq t},  }, \qquad \bm{\theta}_v^t := \qty{ m_v^t, \qty{q_v^{st}, \chi_v^{st}}_{s \leq t} }, \\
  \hat{\bm{\theta}}_u^t &:=  \qty{ \hat{m}_u^t,  \hat{R}^t, \qty{\hat{q}_u^{st}, \hat{\chi}_u^{st}}_{s \leq t},  } , \qquad \hat{\bm{\theta}}_v^t := \qty{ \hat{m}_v^t,\qty{\hat{q}_v^{st},  \hat{\chi}_v^{st}}_{s \leq t} }
\end{align}
and its accumulation as 
\begin{equation}
  \bm{\Theta}_u^t := \bigcup_{s = 1}^t \bm{\theta}_u^s, \quad \bm{\Theta}_v^t := \bigcup_{s = 1}^t \bm{\theta}_v^s, \quad \hat{\bm{\Theta}}_u^t := \bigcup_{s = 1}^t \hat{\bm{\theta}}_u^s, \quad \hat{\bm{\Theta}}_v^t := \bigcup_{s = 1}^t \hat{\bm{\theta}}_v^s.
\end{equation}
Note that $\bm{\Theta}_{\rm RS} \cup \hat{\bm{\Theta}}_{\rm RS} = \bm{\Theta}_u^T \cup \bm{\Theta}_v^T \cup \hat{\bm{\Theta}}_u^T \cup \hat{\bm{\Theta}}_v^T$.
In the successive limit \eqref{eq:successive_limit}, the average generating function $\mathcal{G}(\bm{\Theta}_{\rm RS}, \hat{\bm{\Theta}}_{\rm RS})$ 
is given asymptotically by 
\begin{equation}\label{eq:Extremum_problem}
  \Extr_{\bm{\Theta}_{\rm RS}, \hat{\bm{\Theta}}_{\rm RS}} \mathcal{G}(\bm{\Theta}_{\rm RS}, \hat{\bm{\Theta}}_{\rm RS}) =  \Extr_{\bm{\Theta}_{\rm RS}, \hat{\bm{\Theta}}_{\rm RS}}  \sum_{t = 1}^T \qty[ \beta_u^t \mathcal{G}_u^t(\bm{\Theta}_u^t, \bm{\Theta}_v^t, \hat{\bm{\Theta}}_u^{t}, \hat{\bm{\Theta}}_v^t) + \beta_v^t \mathcal{G}_v^t(\bm{\Theta}_u^{t-1},\bm{\Theta}_v^t, \hat{\bm{\Theta}}_u^{t-1}, \hat{\bm{\Theta}}_v^{t}) ].
\end{equation}
This indicates that the generating function at time $t$, $\mathcal{G}_u^t$, only involves the order parameters up to time $t$, e.g.
the effect of the generating functions $\mathcal{G}_u^s$ at any time $s > t$ and 
$\mathcal{G}_v^s$ at time $s \geq t$ cannot propagate to $\mathcal{G}_u^t$. 
This is a direct consequence of causality in the process defined by \eqref{eq:stochastic_process}. 
The same argument holds for $\mathcal{G}_v^t$, 
which only involves the $u$-indexed order parameters up to time iteration $t$, and the $v$-indexed order parameters up to time iteration $t-1$. 
The subtle difference in time indices between the arguments held by $\mathcal{G}_u^t$ and $\mathcal{G}_v^t$
merely results from the ordering of the alternating procedure in AM \eqref{eq:AM}, 
with the $v$-optimization being followed by the $u$-optimization within a single time index $t$. 

Moreover, 
under the successive limit \eqref{eq:successive_limit}, 
the dominant contribution of $\bm{\theta}_u^t$, 
arises only from the term $ \beta_u^t \mathcal{G}_u^t$, and similarly for $\bm{\theta}_v^t$, only from $ \beta_v^t \mathcal{G}_v^t$, 
since their coefficients $\beta_u^t, \beta_v^t$ are overwhelmingly large compared to $\beta_u^{s}, \beta_v^{s}$ for $s > t$. 
Therefore, the order parameters at the extremum, ${\bm{\theta}}_u^{t,\sharp}, \hat{\bm{\theta}}_u^{t,\sharp},{\bm{\theta}}_v^{t,\sharp }$ and $\hat{\bm{\theta}}_v^{t,\sharp}$ (and their accumulations, $ {\bm{\Theta}}_u^{t,\sharp}, \hat{\bm{\Theta}}_u^{t,\sharp}, { \bm{\Theta}}_v^{t,\sharp}, \hat{\bm{\Theta}}_v^{t,\sharp}$)  
are not determined at once but in a successive manner, each being dependent on the solution of previous iterations: 
\begin{align}
  {\bm{\theta}}_v^{t,\sharp}, \hat{\bm{\theta}}_v^{t,\sharp} &= \arg \Extr_{\bm{\theta}_v^t, \hat{\bm{\theta}}_v^t } \ \mathcal{G}_v^t ( \bm{\theta}_v^t, \hat{\bm{\theta}}_v^t \big| { \bm{\Theta}}_u^{t-1,\sharp}, \hat{\bm{\Theta}}_u^{t-1, \sharp }, { \bm{\Theta}}_v^{t-1, \sharp }, \hat{\bm{\Theta}}_v^{t-1, \sharp }), \\
   \bm{\theta}_u^{t,\sharp}, \hat{\bm{\theta}}_u^{t,\sharp} &= \arg \Extr_{\bm{\theta}_u^{t}, \hat{\bm{\theta}}_u^t } \ \mathcal{G}_u^t ( \bm{\theta}_u^{t}, \hat{\bm{\theta}}_u^t \big| { \bm{\Theta}}_u^{t-1,\sharp}, \hat{\bm{\Theta}}_u^{t-1, \sharp }, { \bm{\Theta}}_v^{t, \sharp }, \hat{\bm{\Theta}}_v^{t, \sharp }).
\end{align}
Therefrom, the functions $\mathcal{G}_u^t$ and $\mathcal{G}_v^t$ are further expressed as 
\begin{align}
  \begin{split}
    &\mathcal{G}_v^t ( \bm{\theta}_v^t, \hat{\bm{\theta}}_v^t \big| { \bm{\Theta}}_u^{t-1,\sharp}, \hat{\bm{\Theta}}_u^{t-1, \sharp}, { \bm{\Theta}}_v^{t-1, \sharp}, \hat{\bm{\Theta}}_v^{t-1, \sharp}) = \frac{q_v^{tt} \hat{q}_v^{tt} - \chi_v^{tt} \hat{\chi}_v^{tt} }{2} - m_v^t \hat{m}_v^t\\
    & \qquad \qquad \qquad \qquad  - \sum_{s < t} ( q_v^{st} \hat{q}_v^{st} + \chi_v^{st} \hat{\chi}_v^{st} ) + \mathcal{S}_v^t(\hat{\bm{\theta}}_v^t \big| \hat{\bm{\Theta}}_v^{t-1, \sharp}) - \kappa \mathcal{E}_v^t\qty({\bm{\theta}}_v^t \big| \bm{\Theta}_u^{t-1, \sharp}, \bm{\Theta}_v^{t-1, \sharp}),
  \end{split}
\end{align}
and 
\begin{align}
  \begin{split}
    &\mathcal{G}_u^t ( \bm{\theta}_u^{t}, \hat{\bm{\theta}}_u^t \big| { \bm{\Theta}}_u^{t-1,\sharp}, \hat{\bm{\Theta}}_u^{t-1, \sharp}, { \bm{\Theta}}_v^{t,\sharp}, \hat{\bm{\Theta}}_v^{t,\sharp}) = \frac{q_u^{tt} \hat{q}_u^{tt} - \chi_u^{tt} \hat{\chi}_u^{tt} }{2} - m_u^t \hat{m}_u^t - R^t \hat{R}^t \\
    & \qquad \qquad \qquad \qquad - \sum_{s < t} ( q_u^{st} \hat{q}_u^{st} + \chi_u^{st} \hat{\chi}_u^{st} ) + \mathcal{S}_u^t(\hat{\bm{\theta}}_u^t  \big| \hat{\bm{\Theta}}_u^{t-1,\sharp} ) - \kappa \mathcal{E}_u^t\qty({\bm{\theta}}_u^t \big| \bm{\Theta}_u^{t-1, \sharp},\bm{\Theta}_v^{t, \sharp}  ). 
  \end{split}
\end{align}
The explicit expressions for $\mathcal{S}_u^t, \mathcal{S}_v^t, \mathcal{E}_u^t,$ and $ \mathcal{E}_v^t$ 
are rather involved, which we provide in the following paragraphs. 

\paragraph{Expression for $\mathcal{S}_v^t, \mathcal{S}_u^t$.} 
The entropic terms $\mathcal{S}_v^t, \mathcal{S}_u^t$ are expressed via Gaussian processes $\qty{\vsf^s}_{s=1}^t$ and $\qty{\usf^s}_{s=1}^t$ respectively, both 
being defined by recursion 
\begin{align}\label{eq:gaussian_process}
  \begin{split}
      \vsf^t &:= \frac{1}{\hat{q}_v^{tt} + \lambda}\Big( x^t_v + \hat{m}_v^t \vsf^\star + \sum_{ t^\prime =1}^{t-1} \hat{q}_v^{t^\prime t} {\vsf}^{t^\prime} \Big),\\
      \usf^t &:= \frac{1}{\hat{q}_u^{tt} + \lambda}\Big( x^t_u + \hat{m}_u^t \usf^\star +  \hat{R}^t \usf^0 + \sum_{ t^\prime = 1}^{t-1} \hat{q}_u^{t^\prime t} {\usf}^{t^\prime} \Big),
  \end{split}
\end{align}
Here, $\vsf^\star, \usf^\star, \usf^0$ are Gaussian random variables given by 
$\usf^\star, \vsf^\star \sim \mathcal{N}(0,1)$, $\usf^0 \sim \mathcal{N}(m_0\usf^\star | 1 - m_0^2)$, 
and $\{x_u^t\}, \{x_v^t\}$ are two independent, centered multivariate Gaussian random variables with covariances $\EE [x_u^t x_u^{t^\prime}] = \hat{\chi}_u^{tt^\prime}$ and 
$\EE [x_v^t x_v^{t^\prime} ]= \hat{\chi}_v^{tt^\prime}$ (assuming $t \leq t^\prime$). The terms of interest $\mathcal{S}_v^t$ and $\mathcal{S}_u^t$ are then given 
simply by 
\begin{equation}
  \mathcal{S}_v^t = \frac{\hat{q}_v^{tt} + \lambda}{2} \EE\qty[ (\vsf^t)^2], \qquad \mathcal{S}_u^t = \frac{\hat{q}_u^{tt} + \lambda}{2} \EE \qty[ (\usf^t)^2].
\end{equation}
Note that the two Gaussian processes $\{\vsf^t\}$ and $\{\usf^t\}$ are independent of each other given a set of order parameters $\bm{\Theta}_{\rm RS}$. 
\paragraph{Expression for $\mathcal{E}_v^t$, $\mathcal{E}_u^t$.}
The energetic terms $\mathcal{E}_v^t, \mathcal{E}_u^t$ are expressed via a sequence of random optimization problems defined by
\begin{align}
  \begin{split}
    L_v^t(w | \{z^s\}_{s < t}, \{w^s\}_{s<t} , h^{t-1}, k^t) &= \frac{w^2}{2 \chi_v^{tt}} + \ell \qty(\phi_u^{t-1} + z^{t-1} ,\phi_v^t + w;y ) , \\
    L_u^t(z | \{z^s\}_{s < t}, \{w^s\}_{s\leq t}, h^t, k^t) &= \frac{z^2}{2 \chi_u^{tt}} + \ell \qty( \phi_u^t + z, \phi_v^t + w^t;y) , 
  \end{split}
\end{align}
where $\phi_u^t(\{z^s\}_{s < t}, h^{t} )$ and $\phi_v^t(\{w^s\}_{s < t}, k^t)$ are given by 
\begin{align}\label{eq:phi_lag}
  \phi_u^t(\{z^s\}_{s < t}, h^{t} ) = h^t + \sum_{s = 1}^{t-1} \frac{\chi_u^{st}}{\chi_u^{ss}} z^s, \quad \phi_v^t(\{w^s\}_{s < t}, k^t) = k^t + \sum_{s = 1}^{t-1} \frac{\chi_v^{st}}{\chi_v^{ss}} w^s.
\end{align}
and $z^t, w^t$ is defined by the recursive relation 
\begin{align}\label{eq:effecive_solution}
  z^t &= \argmin_{z} L_u^t(z | \{z^s\}_{s < t}, \{w^s\}_{s\leq t}, h^t, k^t), \\
   w^t &= \argmin_{w} L_v^t(w | \{z^s\}_{s < t}, \{w^s\}_{s<t} , h^{t-1}, k^t).
\end{align}
The random fields $(k^\star, \bm{k})\in \mathbb{R}^{t + 1}$ and $(h^\star, h^0, \bm{h})\in \mathbb{R}^{t + 2}$ are centered multivariate Gaussian random variables with covariances
\begin{equation}\label{eq:covariances}
  \mqty( 1 & \bm{m}_v^\ten \\
  \bm{m}_v & \vb{Q}_v ) \quad \text{and} \quad
  \mqty( 1 & m_0 & \bm{m}_u^\ten \\
  m_0 &  1 & \bm{R}^\ten \\
  \bm{m}_u & \bm{R} &  \vb{Q}_u ) ,
\end{equation} 
where the vectors $\bm{m}_{u,v}$ and $\bm{R}$ are the concatenation of $m_{u,v}^t$ and $ R^t$ respectively, while 
$\vb{Q}_{u,v}$ is a symmetric matrix with entries $[\vb{Q}_{u,v}]_{st} = q_{u,v}^{\min(s,t), \max(s,t)}$.
The energetic terms $\mathcal{E}_v^t$ and $\mathcal{E}_u^t$ are finally given by the expectation of $L_u^t$ and $L_v^t$ over the random fields: 
\begin{align}
  \mathcal{E}_v^t = \EE \qty[ L_v^t(w^t | \{z^s\}_{s < t}, \{w^s\}_{s<t} , h^{t-1}, k^t)], \quad  \mathcal{E}_u^t = \EE \qty[ L_u^t(z^t | \{z^s\}_{s < t}, \{w^s\}_{s\leq t}, h^t, k^t)]. 
\end{align}

The extremum conditions for $\mathcal{G}_v^t$ are given by the following set of saddle point equations: 
\begin{subequations}\label{eq:u_saddlepoint}
\begin{align}
  %\begin{split}
    m_v^{t} &= \EE\qty[ \vsf^t \vsf^\star], \\
    q_v^{st} &= \EE \qty[ \vsf^s \vsf^t] \quad (s \leq t), \\
    \chi_v^{st} &= \frac{1}{\hat{q}_v^{tt} + \lambda } \qty(\delta_{st} +  \sum_{ \tpr = s}^{t-1} \hat{q}_v^{\tpr t} {\chi}_v^{s \tpr}) \quad (s \leq t), \\
    \hat{m}_v^{t} &= - \kappa \EE \qty[ \frac{\dd^2}{ \dd k^t \dd k^\star} L_v^t ] ,\\
    \hat{q}_v^{st} &=  {\color{black}(2\delta_{st}-1 ) }\kappa \EE \qty[ \frac{\dd^2}{\dd k^s \dd k^t} L_v^t] \quad (s < t) ,\\
    \hat{\chi}_v^{st} &= - \kappa \EE \qty[ \frac{w^t}{\chi_v^{ss}} \partial_2 \ell\qty( \phi_u^{t-1} + z^{t-1}, \phi_v^t + w^t; y ) ] \quad (s < t) ,\\
    \hat{\chi}_v^{tt} &= \kappa \EE \qty[ \qty( \frac{w^t}{\chi_v^{tt}})^2 ].
  %\end{split}
\end{align}
\end{subequations}
Here, $\partial_i \ell$ denotes the partial derivative of $\ell$ with respect to its $i (= 1,2)$-th argument. On the other hand, 
the extremum conditions for $\mathcal{G}_u^t$ are given by the following set of saddle point equations:
\begin{subequations}\label{eq:v_saddlepoint}
  \begin{align}
  %\begin{split}
  \label{eq:m_u}
    m_u^t &= \EE \qty[ \usf^t \usf^\star], \\
    \label{eq:R}
    R^t &= \EE \qty[ \usf^0 \usf^t], \\
    \label{eq:q_u}
    q_u^{st} &= \EE \qty[ \usf^s \usf^t] \quad (s \leq t), \\
    \label{eq:chi_u}
    \chi_u^{st} &= \frac{1}{\hat{q}_u^{tt} + \lambda} \qty(\delta_{st} +  \sum_{\tpr = s}^{t-1} \hat{q}_u^{\tpr t} {\chi}_u^{s \tpr}) \quad (s \leq t), \\
    \hat{m}_u^t &= - \kappa \EE \qty[ \frac{\dd^2}{\dd h^t \dd h^\star} L_u^t], \\
    \hat{R}_u^t &= - \kappa \EE \qty[ \frac{\dd^2}{\dd h^t \dd k^0} L_u^t], \\
    \hat{q}_u^{st} &={\color{black} (2\delta_{st}-1 ) }\kappa \EE \qty[ \frac{\dd^2}{\dd k^s \dd k^t} L_u^t] \quad (s \leq t), \\
    \hat{\chi}_u^{st} &= - \kappa \EE \qty[ \frac{z^t}{\chi_u^{ss}} \partial_2 \ell\qty( \phi_u^t + z^t, \phi_v^t + w^t; y ) ] \quad (s < t), \\
    \hat{\chi}_u^{tt} &= \kappa \EE \qty[ \qty( \frac{z^t}{\chi_u^{tt}})^2 ].
  %\end{split}
\end{align}
\end{subequations}
\rev{
Note that the average over the Gaussian processes $\{\usf^t\}$ can be performed to yield explicit formulae 
for the order parameters $\{m_u^t, R^t, q_u^{st}, \chi_u^{st}\}$. 
The corresponding expressions for \eqref{eq:m_u}, \eqref{eq:R}, \eqref{eq:q_u} and \eqref{eq:chi_u} are given in 
\eqref{eq:m_u_explicit}, \eqref{eq:R_explicit}, \eqref{eq:q_u_explicit} and \eqref{eq:chi_u_explicit} respectively in Appendix \ref{appendix:saddle}, 
with additional details on the formulae for the $v-$order parameters. 
}
\paragraph{Relation to the online setup.}
Our analysis is naturally extendable to the already-known online setup given by \eqref{eq:online_v} and \eqref{eq:online_u},
 in which case the random fields \eqref{eq:random_fields} consists of the inner product between the regressors and covariates given 
at the corresponding time iteration. However, due to the lack of time correlation between the covariates, the random fields are 
effectively expressed only by the order parameters $\{ m^0, m_{u,a}^t, m_{v,a}^t, R^t_a, Q_{u,ab}^{tt}, Q_{v, ab}^{tt} \}$. 
Due to the diagonal nature of the order parameters given as matrices, 
the analysis is much simpler than the full-batch setup. 
In fact, the stochastic process $\{\usf^t, \vsf^t\}$ will lose the memory term as well as off-diagonal correlation in the effective noise $\bm{x}_{u,v}$.

\section{Characterization of the dynamics of alternating minimization}
\label{section:Characterization}

From the above expression of the average generating function, one can obtain a convenient expression for the average of quantities involving the regressors at each iteration. 
This can be done by incorporating the terms in \eqref{eq:AM_observable} into the replica computation done to calculate $\mathcal{G}$. 
For a function $f: \mathbb{R}^T \times \mathbb{R}^T \times \mathbb{R} \times \mathbb{R} \to \mathbb{R}$ acting elementwise on 
$\{\hat{\bm{u}}^t , \hat{\bm{v}}^t\}, \bm{u}^0 , \bm{u}^\star, \bm{v}^\star$, the expectation of $f$ over the data and trajectory of AM is given by 
\begin{equation}\label{eq:effective_trajectory}
 \lim_{N \to \infty} \frac{1}{N} \sum_{i = 1}^{N} \EE_\mathcal{D} \qty[ \langle f( \{\hat{u}^t_i\}, \{\hat{v}^t_i\}, u_i^0, u_i^\star, v_i^\star  ) \rangle_{{\rm AM} | \mathcal{D}, \bm{u}^0 }  ] = \EE \qty[ f( \{\usf^t\}, \{\vsf^t\}, \usf^0, \usf^\star, \vsf^\star ) ].
\end{equation}
This claim indicates that the joint of each element of the regressors, 
initial points, and target vectors at all iterations $t$, $(\{\hat{u}_i^t\},\{ \hat{v}_i^t\}, u_i^0, u_i^\star, v_i^\star)$, 
is statistically equivalent to the joint of the effective random variables $\qty(\{\usf^t\}, \{\vsf^t\}, \usf^0, \usf^\star, \vsf^\star)$ as a population in the large $N$ limit, 
which corresponds to the effective mean-field description of the AM algorithm in the current setup. 
This expression is the first main result of this paper.
These effective random variables are governed by the 
 stochastic process outlined in \eqref{eq:gaussian_process}, with $\{\hat{q}_u^{st}, \hat{q}_v^{st}\}$ 
and the covariance of $(\bm{x}_u, \bm{x}_v)$, $ \{\hat{\chi}_u^{st}, \chi_v^{st}\}$, manifesting the time correlation embedded in this process, 
with their specific values being provided by the solution of the saddle point equations \eqref{eq:u_saddlepoint} and \eqref{eq:v_saddlepoint}. 
By closely examining these order parameters, one can investigate how the memory terms appear and potentially influence the dynamics of the algorithm.

The general expression \eqref{eq:effective_trajectory} allows one to make implications on the order parameters $\{m_u^t, m_v^t\}$, $\{ R^t \}$ and $\{ q_u^{st}, q_v^{st} \}$ in the saddle point equations \eqref{eq:u_saddlepoint} and \eqref{eq:v_saddlepoint}. 
It follows that 
\begin{align}
  \begin{aligned}
    m_u^t &= \lim_{N \to \infty} \frac{1}{N} \EE_\mathcal{D} [ (\bm{u}^\star)^\ten \hat{\bm{u}}^t],\qquad \qquad   &q_u^{st} = \lim_{N\to \infty} \frac{1}{N} \EE_\mathcal{D}  [ (\hat{\bm{u}}^s)^\ten \hat{\bm{u}}^t ], \\
    m_v^t &= \lim_{N \to \infty} \frac{1}{N} \EE_\mathcal{D}  [ (\bm{v}^\star)^\ten \hat{\bm{v}}^t],\qquad\qquad  &q_v^{st} = \lim_{N\to \infty} \frac{1}{N} \EE_\mathcal{D}  [ ( \hat{\bm{v}}^s)^\ten \hat{\bm{v}}^t ], \\
    R^t &= \lim_{N \to \infty} \frac{1}{N} \EE_\mathcal{D}  [ (\bm{u}^0)^\ten \hat{\bm{u}}^t ].
  \end{aligned}
\end{align}
Therefore, $\{m_u^t, m_v^t\}$ and $\{R^t\}$ expresses the overlap between the regressors at each iteration and the target vectors and initial points respectively (note that only the $\bm{u}$ variable is given an initialization), 
while $\{q_u^{st}, q_v^{st}\}$ expresses the overlap between the regressors at different iterations. 
These macroscopic quantities will be used to characterize the dynamics of the AM algorithm in this study. 

\paragraph{Generic factorized priors on the target vectors.}
\label{subsection:generic_priors}
\rev{
One can consider a more generic factorized distribution for the target vectors:
$u_i^\star \sim p_u(u_i^\star)$ and $v_i^\star \sim p_v(v_i^\star)$, where $p_u, p_v$ are arbitrary distributions. 
While our analysis exclusively focused on the case where $p_u = p_v = \mathcal{N}(0,1)$, 
it should be noted that the generating functional and the saddle point equations is 
valid as long as $p_u$ and $p_v$ are centered and possess a unit variance. 
In fact, the derivation of $\mathcal{G}_{\rm RS}(\bm{\Theta}_{\rm RS}, \hat{\bm{\Theta}}_{\rm RS})$ 
does not utilize the specific form of $p_u$ and $p_v$, but only its second moment. 
However, it should be noted that the mean-field description \eqref{eq:effective_trajectory} will be 
altered by the choice of $p_u$ and $p_v$, as the effective random variables $\usf^\star $ and $\vsf^\star$
must be distributed according to $p_u$ and $p_v$, respectively. 
}

\subsection{Impossible retrieval from random initialization}
\label{subsection:random_init}

Since the parameters are determined in a successive manner, one can use mathematical induction 
on the target vector overlap $m_u^t$ and $m_v^t$ 
 to prove the following for AM initialized in a completely random manner ($m_0 = 0$). 
 \begin{claim}\label{claim:random_init}
  Suppose $m_0 = 0$. Then, $m_u^t = m_v^t = 0$ for finite $t$ and finite $\kappa$. 
\end{claim}
The full proof is given in Appendix \ref{appendix:claim1_proof}. 
This indicates that one must have an initialization point with finite correlation with the target 
in order to retrieve anything under finite number of iterations $t$ and finite $\kappa$. 
\rev{We note that this does not exclude the possibility of retrieval from random initialization 
under $t$ and $\kappa$ diverging with $N$, which is in fact possible under an online setup as shown in \cite{Verchand22}.}

\section{Numerical comparison with finite size experiments}
\label{section:numerical}

\begin{figure}
  \centering
  \includegraphics[width=1.0\textwidth]{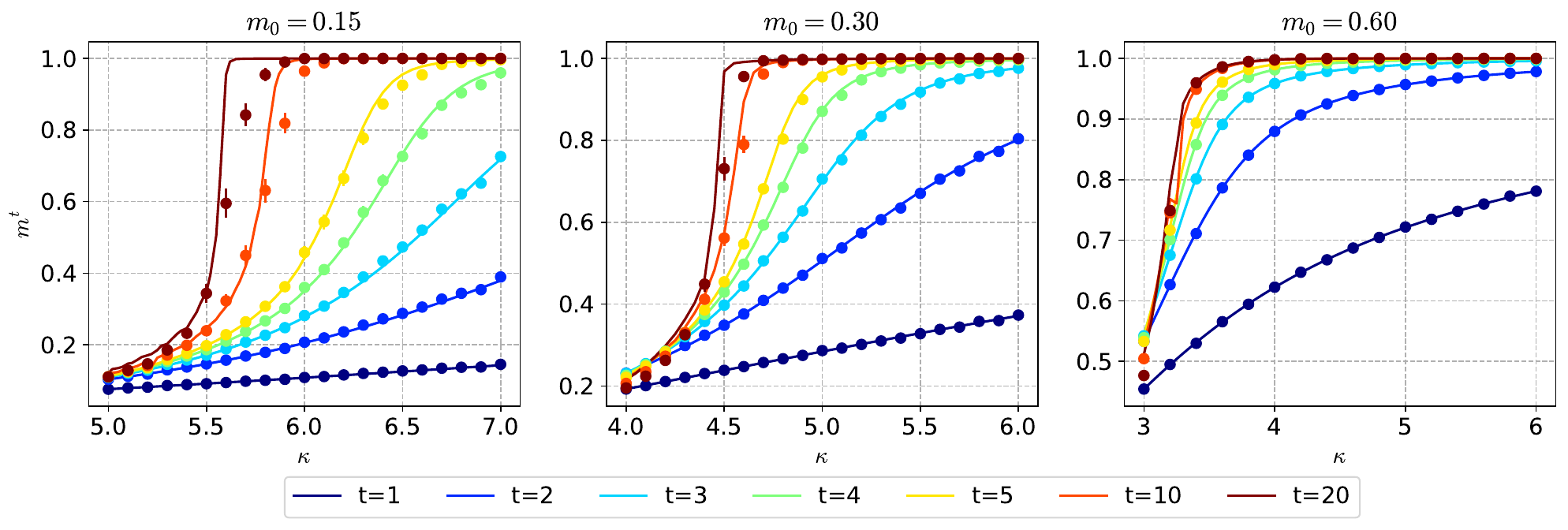}
  \caption{Comparison of the theoretical value (solid line) of $m^t$ and the empirical value (markers) obtained from experiments for $N = 16000$. 
  The theoretical value was obtained by solving the fixed-point equations given in \eqref{eq:u_saddlepoint} and \eqref{eq:v_saddlepoint}.
  The empirical value was obtained by taking the mean over 64 random configurations of $\mathcal{D}$. Error bars represent the 
  standard error of the mean. 
  }
  \label{fig:asymptotic_dynamics}
\end{figure}
In this section, we analyze the behavior of the AM algorithm for bilinear regression by numerically solving the saddle point equations \eqref{eq:u_saddlepoint} and \eqref{eq:v_saddlepoint}. 
Comparisons with experiments on finite size systems are also included. 
 Hereforth, we focus on the quadratic biconvex loss, i.e. $\ell(a,b;y) = \frac{1}{2} (y-ab)^2$, 
 \rev{in which the explicit update procedure of AM in \eqref{eq:AM} is given by 
 \begin{align}
  \bhv^{t} &= ( \mathbf{B}^\ten (\vb{D}_u^{t-1})^2 \vb{B} + \lambda \vb{I}_N )^{-1} \vb{B}^\ten (\vb{D}_u^{t-1})^\ten \bm{y}, \ &&\text{where} \quad \vb{D}^{t-1}_u = {\rm diag}(\vb{A} \bm{u}^{t-1}), \label{eq:AM_v} \\
  \bhu^{t} &= ( \vb{A}^\ten (\vb{D}_v^{t})^2 \vb{A} + \lambda \vb{I}_N )^{-1} \vb{A}^\ten (\vb{D}_v^{t })^\ten \bm{y} \label{eq:AM_u}, &&\text{where} \quad \vb{D}^{t}_v = {\rm diag}(\vb{B} \bm{v}^{t}),
\end{align}
and the matrices $\vb{A}, \vb{B} \in \mathbb{R}^{P\times N}$ are stacked versions of $\{ \bm{A}_\mu, \bm{B}_\mu \}_{\mu = 1}^P$. 
These updates, consisting of basic linear algebra operations, can be performed efficiently using GPUs, in which extensive finite-size experiments can be conducted. }
Also note that a finite $\lambda > 0$ is necessary for the target optimization function \eqref{eq:target_problem} to have a unique minimum, 
as it is invariant under the transformation $(\bm{u}, \bm{v}) \to ( C \bm{u}, \bm{v} / C )$ for any $C$ if $\lambda = 0$. 
To avoid possible complications arising from these degeneracies, we set $\lambda$ to a finite but small value, $\lambda = 0.01$, 
for all experiments. 
We refer the reader to Appendix \ref{appendix:saddle} for a detailed explanation on how to numerically solve the saddle point equations \eqref{eq:u_saddlepoint} and \eqref{eq:v_saddlepoint}.

Recall that the main quantity of interest is the product cosine similarity $m^t$, which we rewrite here for sake of convenience: 
\begin{equation*}
  m^t := \lim_{N\to \infty} \EE_{\Data} \qty[ m_N^t(\Data)],\rev{ \quad \text{where} \quad m_N^t(\Data) = \frac{1}{N} \frac{ ( \hat{\bm{u}}^t )^\ten \bm{u}^\star (\hat{\bm{v}}^t)^\ten \bm{v}^\star }{ \norm{\hat{\bm{u}}^t} \norm{\hat{\bm{v}}^t} }}.
\end{equation*} 
In order to assess this quantity, it is important to note that in the limit $N \to \infty$, the norm of the regressors $\bm{u}^t$ and $\bm{v}^t$ 
given a realization of data $\Data$ typically concentrate on its average. 
This concept of concentration is called \textit{self-averaging}, which has been observed and proven
 in convex optimization \cite{Bayati12,Stojnic13, Thrampoulidis18, Miolane21} and Bayes optimal inference \cite{Jean18,Jean19}. 
Since our problem is essentially a sequence of convex optimization problems, 
we expect the same phenomenon to hold. This observation allows us to 
evaluate $m^t$ as 
\begin{equation}
    m^t = \frac{ m_u^tm_v^t }{\sqrt{q_u^{tt} q_v^{tt}}}.
\end{equation}

\begin{figure}
  \centering
  \includegraphics[width = 1.0\textwidth]{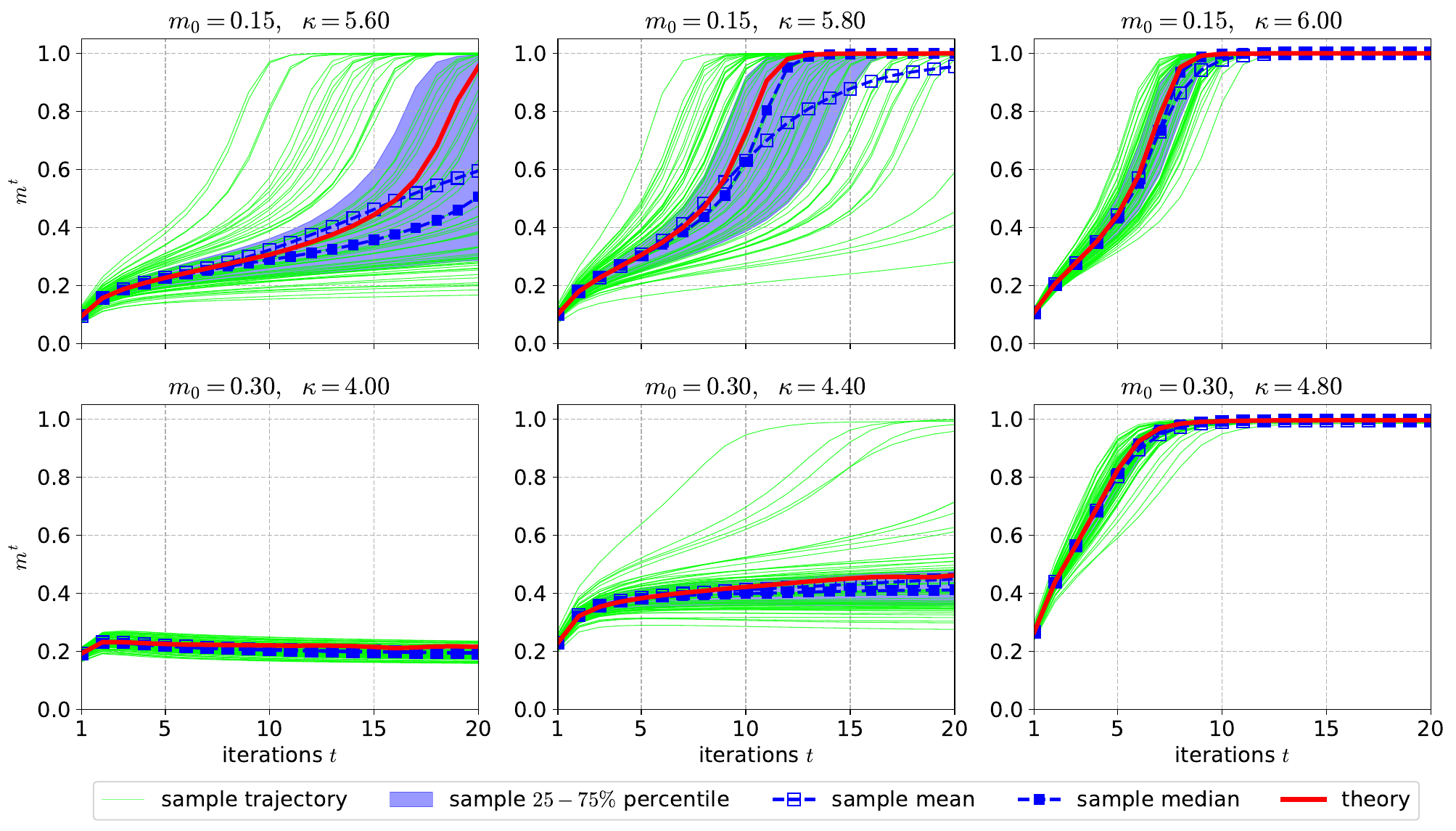}
  \caption{Detailed dynamics of $m^t$ for $m_0 = 0.15$ (top) and $m_0 = 0.30$ (bottom) for various values of $\kappa$. 
  The thin green lines correspond to all 64 independent runs of AM with system size $N = 16000$. 
  We see that the variance of $m^t$ is large for small $\kappa$ and $m_0$, with both mean and median of the population of trajectories 
  deviating from the theoretical value. 
  }
  \label{fig:Trajectories}
\end{figure}

\subsection{Time evolution of the product cosine similarity}

In figure \ref{fig:asymptotic_dynamics}, we compare the value $m^t$ obtained from theory and its empirical counterpart, $\EE_{\Data}[ m^t_N(\Data)]$, with $N = 16000$ for different values of $\kappa$ and $t$. 
Recall that $\kappa$ is the sample complexity; $\kappa=P/N$, where $P$ and $N$ represent the sample size and the dimension of the target vectors, respectively.
The empirical values were obtained by taking the mean over 64 random configurations of $\mathcal{D}$. 
The results from theory and experiment agree well excluding the case when $\kappa$ is small for $m_0 = 0.15$ and $0.30$, which suggests that our effective description basically explains the behavior of AM algorithm correctly. 

To investigate the inconsistency in the case of $m_0=0.15$ and $0.3$ for small $\kappa$ in more detail, 
in figure \ref{fig:Trajectories} we show the detailed dynamics of $m^t$ for $m_0 = 0.15$ and $0.30$ for various values of $\kappa$.
For $\kappa = 5.60$ and $5.80$ for $m_0 = 0.15$, and $\kappa = 4.40$ for $m_0 = 0.30$, 
we see that a typical trajectory of $m^t$ cannot be identified from the experimental values. 
Trajectories from theory and experiment only agree for a small number of iterations, where the variance in the empirical value is small. 
This indicates that even for the system size $N = 16000$, the self-averaging effect is not strong enough for the theoretical value, which 
was derived assuming self-averaging, to be a good approximation of finite-size behavior. 
 Such large finite-size effects, commonly observed when a physical system is close to a critical point, 
 suggests the existence of an \textit{algorithmic} critical point for AM, where the algorithm 
bifuricates into two different dynamical behaviors; one where the algorithm converges to an informative fixed point ($m^t \simeq 1$), 
and the other where the algorithm converges to fixed point with mediocre signal recovery performance. 
In fact, this behavior is already observed in figure \ref{fig:Trajectories} for $m_0 = 0.30$, where $m^t$ seems to converge 
to a small value and $1.0$ for $\kappa = 4.00, 4.60$ respectively, with an indecisive behavior for a value of $\kappa$ in between ($\kappa = 4.40$).
Unfortunately, precisely investigating this critical point (as well as its existence) is not possible with the current analysis, 
as the framework at hand is incapable of handling infinite iteration $t$, and we leave this as a future work.

Nevertheless, this suggests that initialization techniques such as spectral methods \cite{Li19, Charisopoulos20}, may be highly effective. 
A sufficiently accurate initial state, combined with adequate sample complexity within the alleged algorithmically informative phase, 
would drive the algorithm towards the informative fixed point where near-perfect recovery of the target vectors is achieved. 

\subsection{Finite-size effects and algorithmic critial points}

\begin{figure}
  \centering
  \includegraphics[width = 1.0\textwidth]{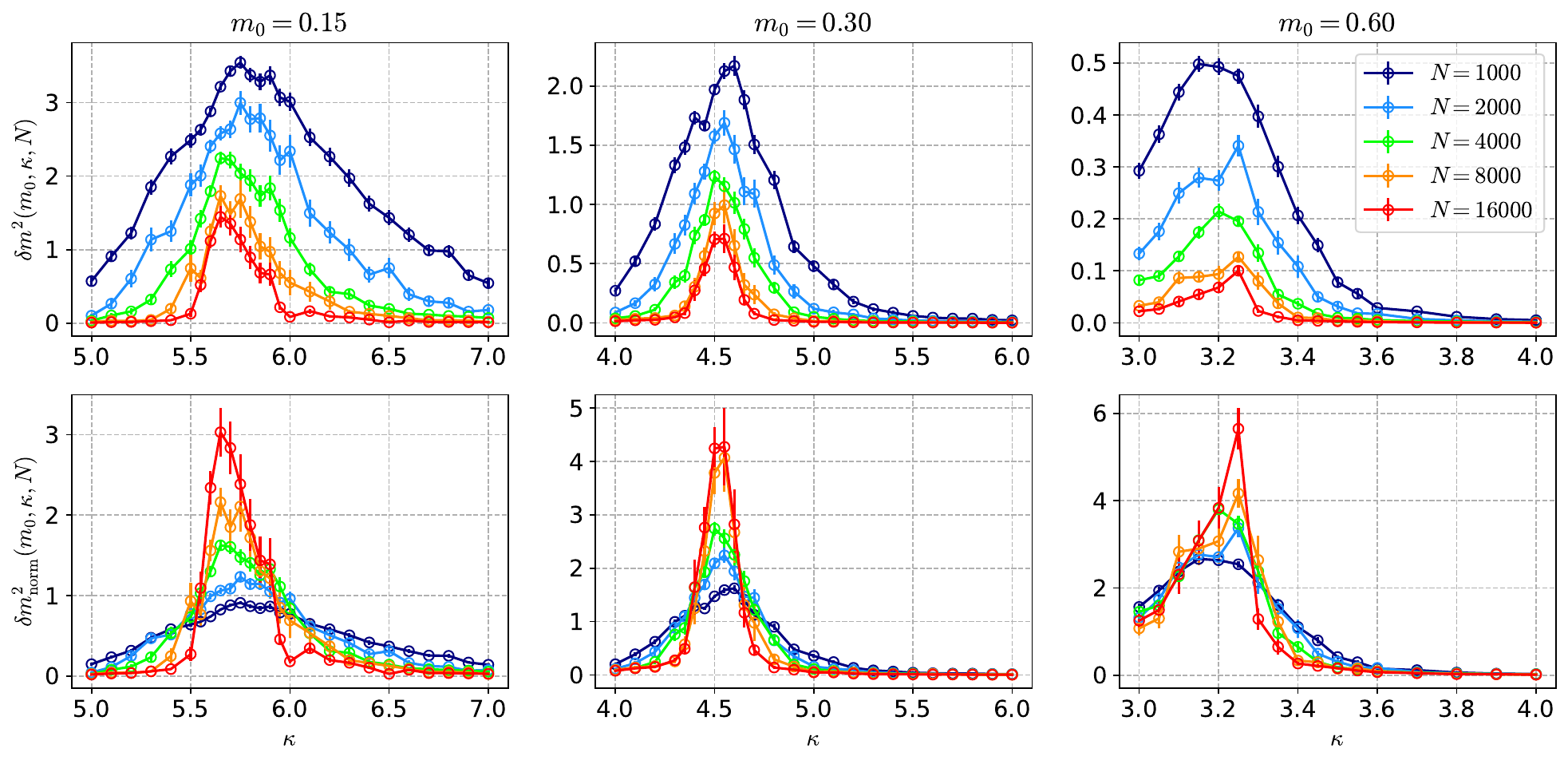}
  \caption{
    \rev{Values of $\delta m^2 (m_0, \kappa, N)$ (upper panel) and its normalized counterpart $\delta_{\rm norm}^2 (m_0, \kappa, N)$ (lower panel) for $m_0 = 0.15$ (left), $0.30$ (middle) and $0.60$ (right) as a function of $\kappa$ for various values of $N$.
  The average over $\Data$ was taken over 1024, 256, 256, 64 and 64 random configurations for $N =$ 1000, 2000, 4000, 8000 and 16000 respectively. 
  Error bars represent the standard error of the mean. }
  }
  \label{fig:Delta_M}
\end{figure}
\rev{
While the analytical framework does not provide a direct way to investigate an algorithmical critical point, 
numerical experiments can provide insights on the existence of such a point. 
Indeed, one can aniticipate that near a critical point, the deviation of the experimental value $m^t_N(\Data)$ from its theoretical counterpart 
$m^t$ should inhibit large finite-size effects, 
as the system is close to a bifurcation point in which the algorithm switches between two different convergence behaviors. 
We therefore consider the squared sum of such deviation summed over $t = 1,\ldots, 20$ as a measure for the finite-size effect: 
\begin{equation}
  \delta m^2 (m_0, \kappa, N) = \EE_{\Data} \qty[ \sum_{t = 1}^{20} ( m^t - m^t_N(\Data) )^2 ].
\end{equation}
The upper panel of Figure \ref{fig:Delta_M} reveals characteristic peaks in $\delta m^2(m_0, \kappa, N)$ as a function of $\kappa$, 
whose positions shift with $m_0$. 
These peaks diminish but also become increasingly sharp for larger $N$, 
suggesting the presence of critical points $\kappa_c(m_0)$ 
where the finite-size effects are maximized. 
Further considering a normalized version of $\delta m^2(m_0, \kappa, N)$ with respect to $\kappa$, i.e. 
\begin{equation}
  \delta m^2_{\rm norm}(m_0, \kappa, N) = \frac{\delta m^2(m_0, \kappa, N)}{ I(m_0, N) }, \quad I(m_0, N) = \int_{\kappa_{\rm min}}^{\kappa_{\rm max}} \dd \kappa\ \delta m^2(m_0, \kappa, N) ,
\end{equation}
with $(\kappa_{\rm min}, \kappa_{\rm max}) = (5.0, 7.0), (4.0, 6.0)$ and $(3.0, 4.0)$ for $m_0 = 0.15, 0.30$ and $0.60$ respectively, 
reveals a peak structure whose height increases with $N$ (Figure \ref{fig:Delta_M}, lower panel).
This behavior strongly suggests the existence of algorithmic critical points where the dynamics of AM transitions between different convergence behaviors. 
}

\begin{figure}
  \includegraphics[width = 0.8\textwidth]{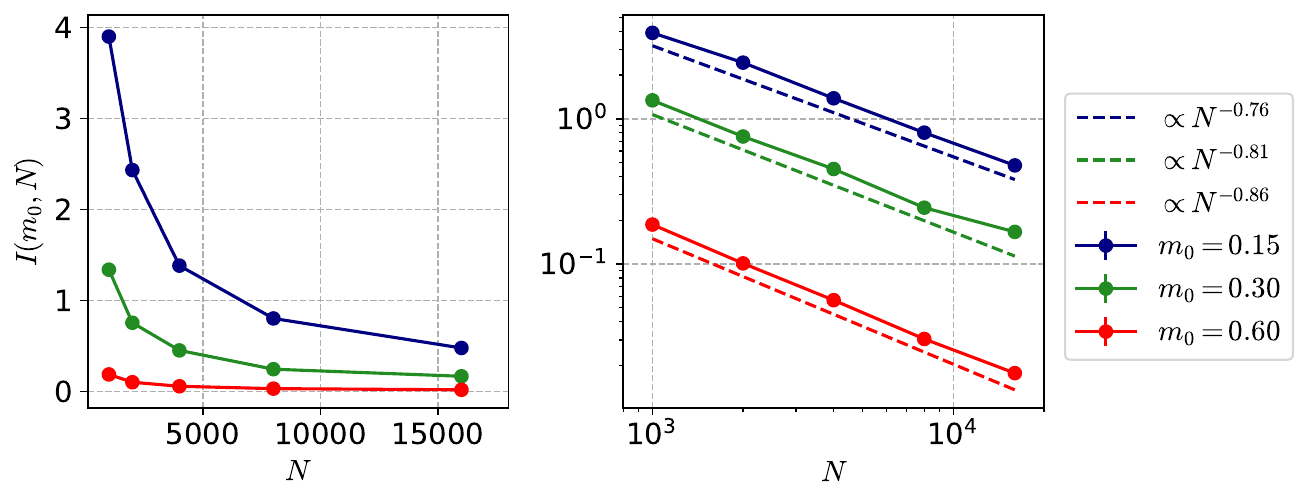}
  \caption{\rev{Integral of $\delta m^2(m_0, \kappa, N)$ over $\kappa$ for $m_0 = 0.15, 0.30$ and $0.60$ as a function of $N$ 
  in normal scale (left) and log-log scale (right). 
  The integral, calculated using the trapezoidal rule, was taken over the region displayed in figure \ref{fig:Delta_M}. 
  Error bars represent the standard error of the mean, which are too small to be visible.}
  }
  \label{fig:Delta_M_scaling}
\end{figure}

\rev{
To investigate the overall effect of $m_0$ on finite-size effects, we plot 
the normalization constant $I(m_0, N)$ over $\kappa$ for $m_0 = 0.15, 0.30$ and $0.60$ in Figure \ref{fig:Delta_M_scaling}. 
While this normalization constant decreases with increasing $m_0$ and $N$, 
the power decay of the integral with respect to $N$ appear to be universal across different values of $m_0$, 
taking the form of $N^{-\alpha}$ with $\alpha \simeq 0.8$, suggesting that $m_0$ only affects the overall magnitude of the finite-size effects, but its scaling behavior. 
}

 \subsection{Time correlation of the dynamics}
\begin{figure*}[!h]
    \centering 
    \begin{minipage}[t]{0.47\textwidth}
        \centering
        \includegraphics[width=1.0\textwidth]{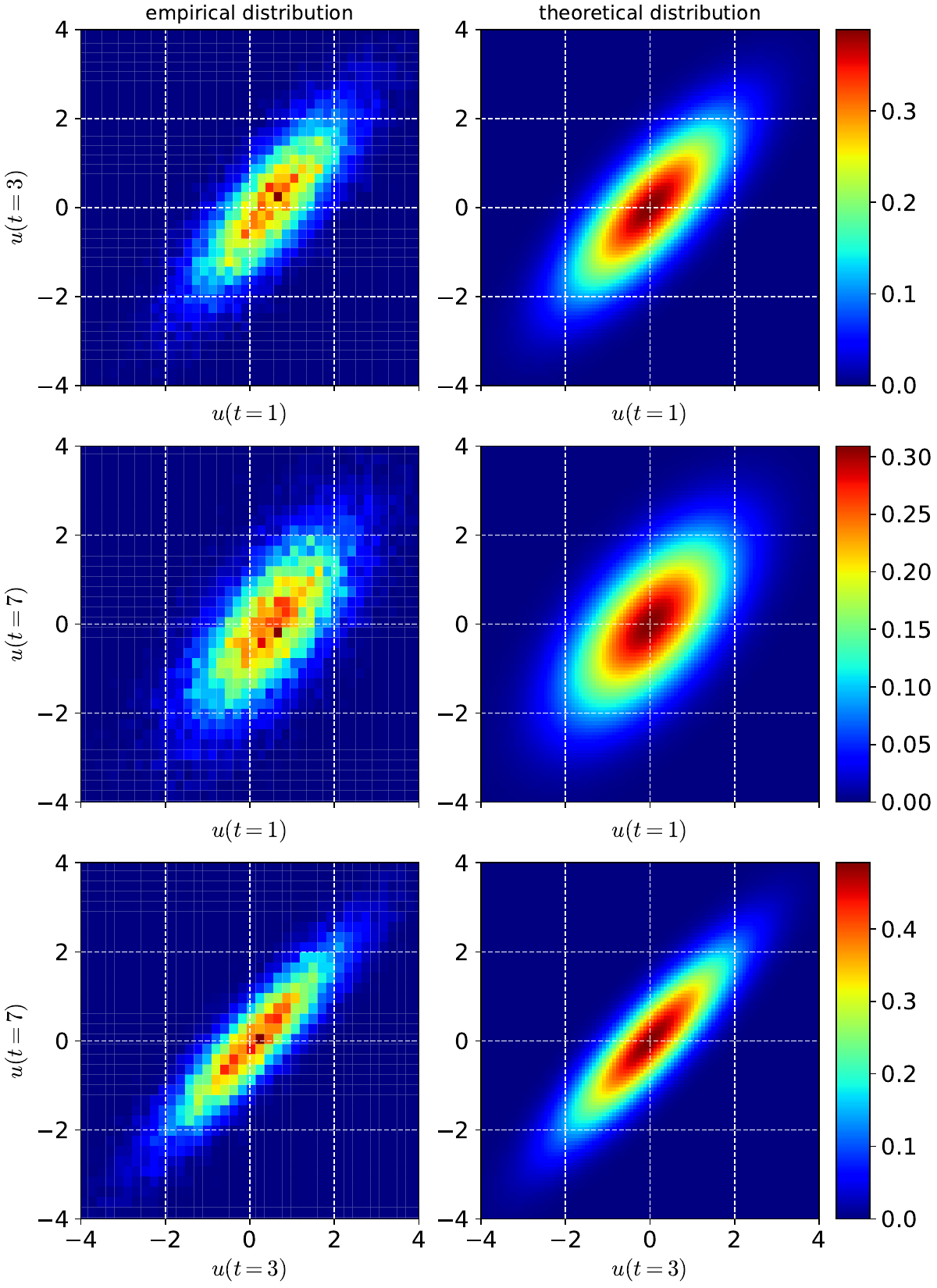}
    \end{minipage}\hfill
    \begin{minipage}[t]{0.47\textwidth}
        \centering
        \includegraphics[width=1.0\textwidth]{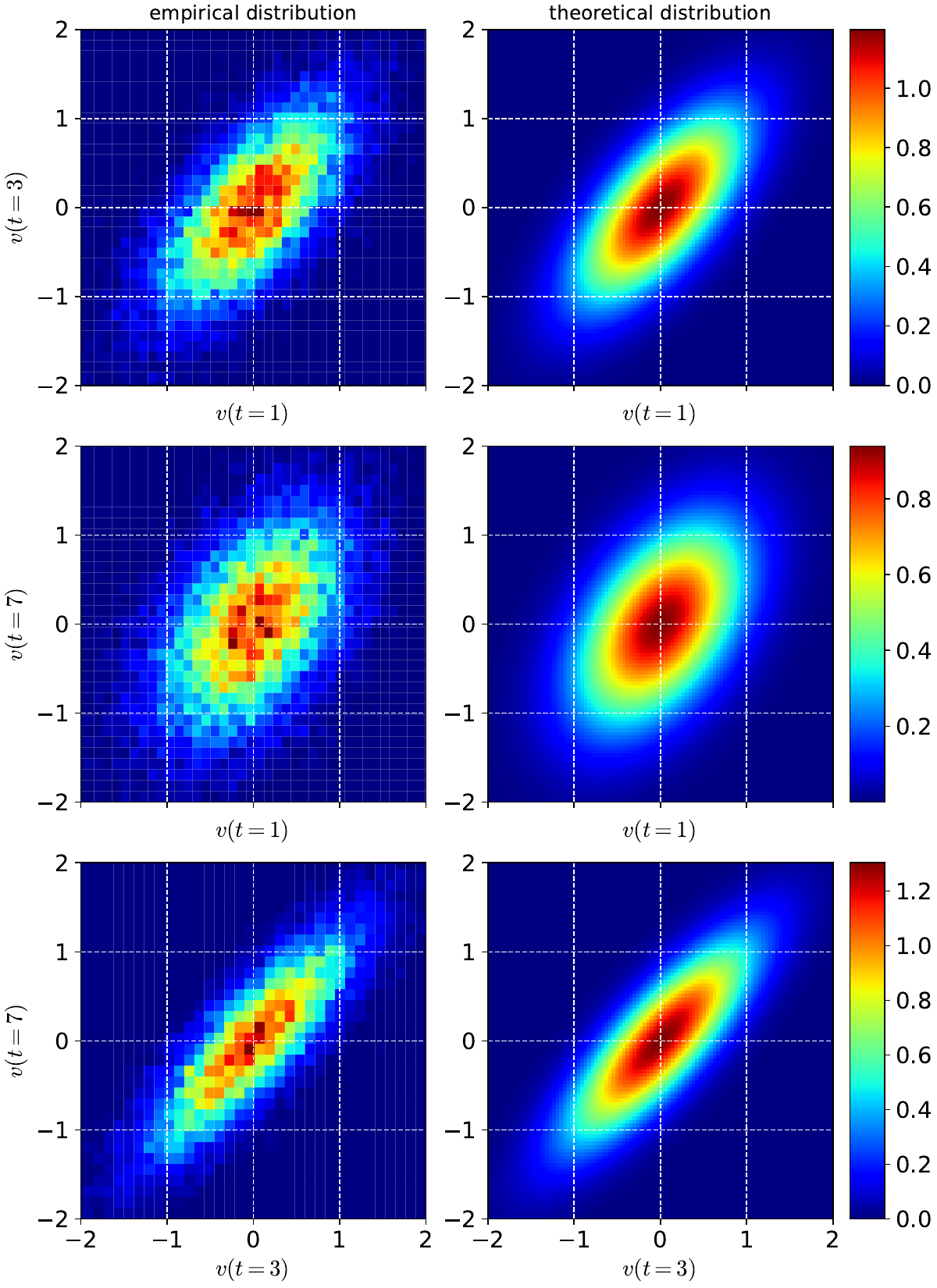}
    \end{minipage}
    \caption{Comparison of the empirical distribution of $\bm{u}^t$ and its theoretical counterpart $\usf^t$ (left), 
    and the empirical distribution of $\bm{v}^t$ and $\vsf^t$ (right) for $t = 1, 3, 7$ and $m_0 = 0.30, \kappa = 5.0$. 
    The empirical distribution was obtained from a single random instance of size $N= 16000$. }
    \label{fig:distribution}
\end{figure*}

The empirical distribution of $\bm{u}^t$ for a single random instance is compared with its theoretical counterpart $\usf^t$. 
In Figure \ref{fig:distribution} we show the joint distribution of $\bm{u}^t$ and $\usf^t$ for $t = 1, 3, 7$ for $m_0 = 0.30, \kappa = 5.0$ and $N = 16000$ for experiments. 
Note that from the underlining Gaussian process describing the asymptotic dynamics, the joint distribution of $(\usf^s, \usf^t)$ is given by a multivariate Gaussian distribution with zero mean and 
covariance 
\begin{equation}
  \mqty( q_u^{ss} & q_u^{st} \\ q_u^{st} & q_u^{tt} ).
\end{equation}
As evident from the plot, the empirical distribution of $\bm{u}^t$ is in good agreement with the theoretical distribution of $\usf^t$, 
even for a single random instance. 

\begin{figure*}[!h]
  \centering 
  \includegraphics[width = 1.0\textwidth]{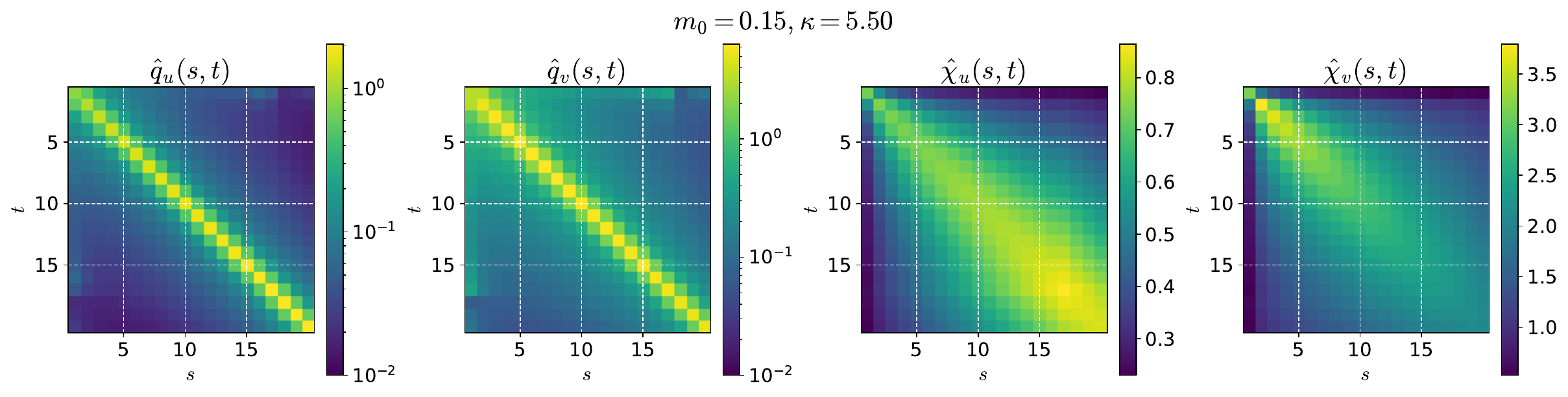}
  \includegraphics[width = 1.0\textwidth]{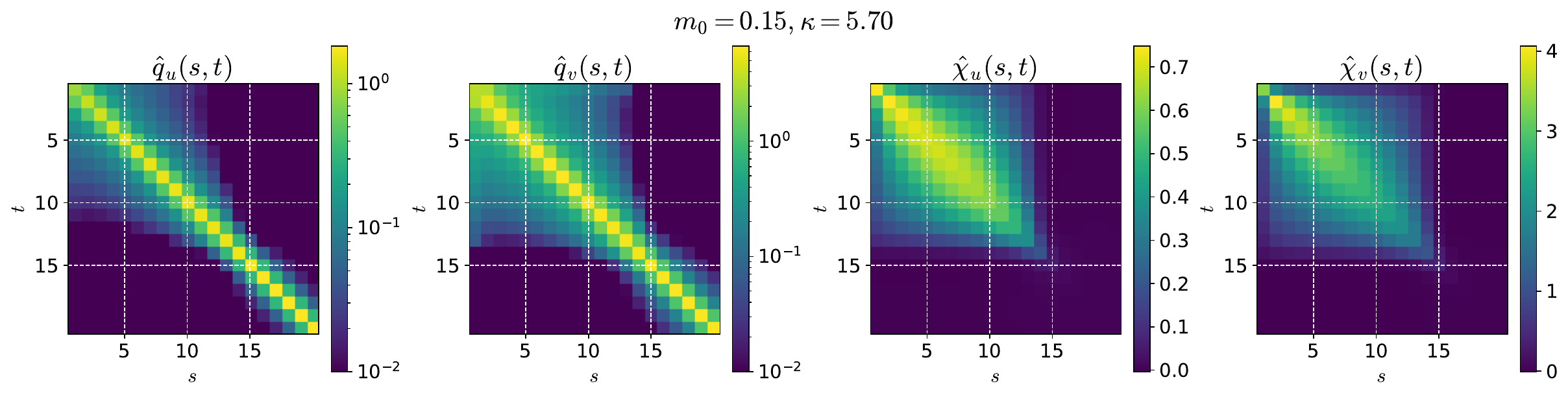}
  \includegraphics[width = 1.0\textwidth]{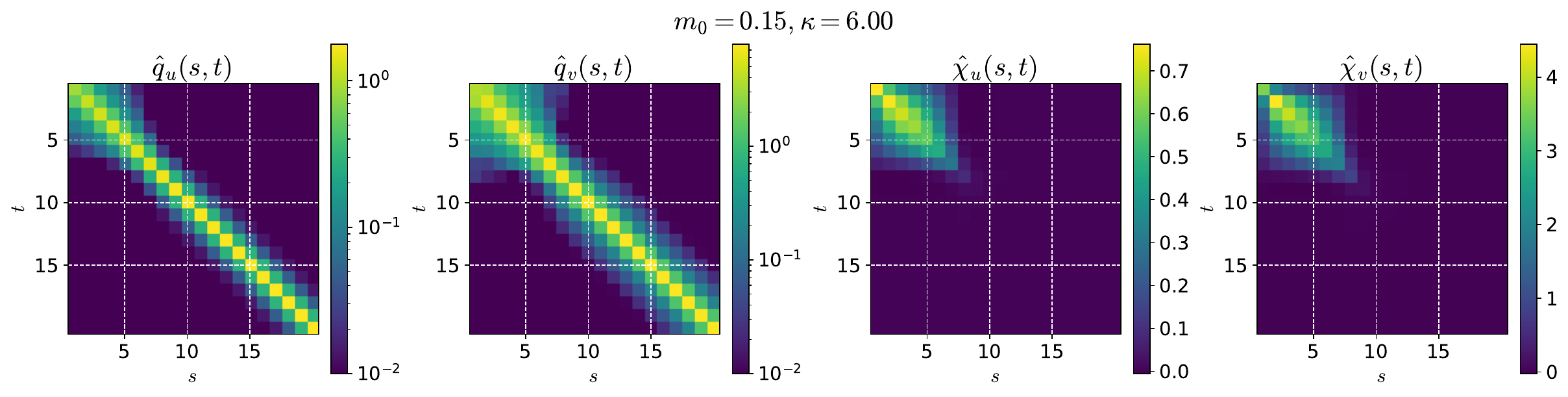}
  \caption{
  Correlation matrix of $\hat{q}_u^{st}, \hat{q}_v^{st}, \hat{\chi}_u^{st}, \hat{\chi}_v^{st}$ for $m_0 = 0.15$ and $\kappa = 5.50, 5.70$ and $6.00$.
  Note that $q_u^{st}$ and $q_v^{st}$, only defined for $s \leq t$, is symmetrized, and has logscale colorbars for sake of visualization. }
  \label{fig:time_matrix}
\end{figure*}

In Figure \ref{fig:time_matrix} we plot the matrices $\{\hat{q}_u^{st}\}, \{\hat{\chi}_u^{st}\}, \{\hat{q}_v^{st}\}$ and $\{\hat{\chi}_v^{st}\}$,
which induce the memory effect in the Gaussian process \eqref{eq:gaussian_process}. 
Note that although $\{\hat{q}_u^{st}, \hat{q}_v^{st}\}$ are only defined for $s < t$, we symmetrize the matrices for sake of visualization. 
For cases where retrieval of the signal is possible within the span of $t \leq T = 20 \ (\kappa = 6.0, 5.70)$, 
the external noise term $\{x_u^t\}$, whose covariance is given by $\{\hat{\chi}_u^{st}\}$, 
quickly disappears from a certain iteration. 
The matrix $\{\hat{q}_u^{st}\}$, which resembles the 
lag term in the Gaussian process, 
also holds a significant time correlation during the retrieval process. 
However, even after the signal has been recovered and the external noise term has disappeared, 
$\{\hat{q}_u^{st}\}$ still holds a short-term memory effect. 
The same behavior is qualitatively observed for the $v$-matrices, but they are not quantitatively identical; 
this is because the AM algorithm is not symmetric in the sense that only $\bm{u}$ is given an initialization, and 
$\bm{v}^t$ is always calculated ahead of $\bm{u}^t$. 
Nevertheless, this suggests that in the earlier iterations of the algorithm, strong memory effects 
appear both in the form of time-correlated noise and lag which acts as an external force driving the dynamics to a fixed point. 
After recovery has been achieved, 
 the stationary Gaussian process has no external noise $\bm{x}_{u,v}$ but still possesses a short-termed memory effect. 

\section{Conclusion and discussion}
In this work, we have obtained a closed-form expression for the asymptotic dynamics of AM using the replica method. 
Our result conjectures that the regressors at each iteration can be statistically characterized by a stochastic process, 
shedding light on the algorithm's effective memory dependency.
Numerical results suggest that our analysis captures the asymptotic dynamics. 
Moreover, examination of memory terms in the stochastic process reveals that only short-term memory influence dynamics at later iterations, 
in contrast to a more pronounced long-term memory dependency during its early evolution.

From a technical viewpoint, our analysis can be extended to other types of loss functions and iterative algorithms under random data, 
opening exciting directions for future exploration. 

Moreover, while the initialization setup given in \eqref{eq:initialization} is purely conventional, we stress that 
the analysis in our work is extendable to more realistic spectral initializations \cite{Chi19,Xiaodong19,Charisopoulos20}. 
In general, this would only require a modification of the distribution $P(\bm{u}_0 | \bm{u}_\star)$ to $P(\bm{u}_0 | \mathcal{D})$, 
which can be handled in the same manner as the current analysis.  
In addition, it will be interesting investigate whether a phase transition from a retrieval to a non-retrieval phase 
exists for the AM algorithm; such algorithmic critical points are known to exist in methods such as approximate message passing in terms of sample complexity \cite{Krzakala12, Maillard20}. 
While we provide some numerical evidence suggesting its presence, 
a more established study would require the analysis of AM in the limit of both $N \to \infty$ and $T \to \infty$, which is a challenging task left for future work.

\section*{Acknowledgements}
The authors would like to thank Yoshiyuki Kabashima for insightful discussions and comments. 
Computational resource of AI Bridging Cloud Infrastructure (ABCI) provided by National Institute of Advanced Industrial Science and Technology (AIST) was used.
% TODO: include author contributions
% TODO: include funding information
\paragraph{Funding information}
This work is supported by JSPS KAKENHI Grant Nos. 22KJ1074 (KO), 23K16960, 21K21310 (TT), 20H00620 (TT, KO), and JST CREST Grant Number JPMJCR1912 (TT, KO). 
\begin{appendix}
\numberwithin{equation}{section}

\section{Derivation of replica symmetric average generating function}
\label{appendix:replica_derivation}
\subsection{Evaluation of the state density term}

We start by evaluating the state density term $\mathcal{V}(\vb{\Theta})$ under the replica symmetric ansatz. 
By using the Fourier representation of the delta function, we have
\begin{align}
  \begin{split}
  \prod_{a, b} &\delta\qty(NQ_{u,ab}^{st} - (\bm{u}_a^s)^\ten \bm{u}_b^t) = \prod_{a = 1}^n \delta\qty( N q_u^{st} - (\bm{u}_a^s)^\ten \bm{u}_a^t) \prod_{a \neq b} \delta\qty( N q_{u}^{st} - \frac{N\chi_u^{st}}{\beta_u^s} - (\bm{u}_a^s)^\ten \bm{u}_b^t)\\
  &\propto \int \dd \hat{q}_u^{st} \dd \hat{\chi}_u^{st} \exp N n \beta_u^{t} \qty[ \qty( 1 - \frac{\delta_{st}}{2} ) ( q_u^{st} \hat{q}_u^{st} + (n-1) \chi_u^{st}\hat{\chi}_u^{st} )   ] \\
  &\times \exp  \qty( 1 - \frac{\delta_{st}}{2} ) \qty[  -\beta_u^{t} \hat{q}_u^{st} \sum_{a =1}^n (\bm{u}_a^s)^\ten \bm{u}_a^t + \beta_u^s \beta_u^t \hat{\chi}_u^{st} \qty( \sum_{a =1}^n \bm{u}_a^{s} )^\ten \qty( \sum_{a =1}^n \bm{u}_a^{t} ) ].
  \end{split}
\end{align}
Taking the product over pairs of $s \leq t \ (\leq T)$ offers
 \begin{align}
  \begin{split}
  &\int \prod_{s \leq t}^T \dd \hat{q}_u^{st} \dd \hat{\chi}_u^{st}  \exp \qty{ N n \sum_{t = 1}^T \beta_u^t \qty( \frac{q_u^{tt} \hat{q}_u^{tt} - \chi_u^{tt}\hat{\chi}_u^{tt}}{2} + \sum_{ s < t} (q_u^{st} \hat{q}_u^{st} - \chi_u^{st}\hat{\chi}_u^{st}) + O(n)) } \\
  \times &\exp \qty{ \sum_{a = 1}^n \sum_{t = 1}^T \beta_u^t  \qty( -\frac{1 }{2} \hat{q}_u^{tt} \qty( \bm{u}_a^t)^\ten \bm{u}_a^t -  \sum_{s < t}^T  \hat{q}_u^{st} \qty( \bm{u}_a^s)^\ten \bm{u}_a^t) + \sum_{s,t = 1}^T \beta_u^s \beta_u^t \chi_u^{st} \qty( \sum_{a =1}^n \bm{u}_a^{s} )^\ten \qty( \sum_{a =1}^n \bm{u}_a^{t} ) }.
  \end{split} 
\end{align}
In order to decouple the last term in the exponential of the above expression with respect to the replica indices, we introduce a multi-dimensional Hubbard Stratonovich transformation, based on the following trivial identity: 
\begin{equation}
  \EE_{\bm{x} \sim \mathcal{N}(\bm{0}, \vb{\Sigma})}\qty[ e^{\bm{a}^\ten \bm{x}} ] = e^{\frac{1}{2} \bm{a}^\ten \vb{\Sigma} \bm{a}}.
\end{equation}
Applying this identity with respect to the time dimension, one obtains 
\begin{equation}
  \exp \qty{  \frac{1}{2}\sum_{s,t = 1}^T \beta_u^s \beta_u^t \chi_u^{st} \qty( \sum_{a =1}^n \bm{u}_a^{s} )^\ten \qty( \sum_{a =1}^n \bm{u}_a^{t} ) } =\prod_{i = 1}^N \EE_{\bm{x}_u \sim \mathcal{N}(\bm{0}, \vb{\chi}_u)} \qty[ \exp \qty{ \sum_{t = 1}^T \beta_u^t x_u^t \sum_{a =1}^n u_{a,i}^t } ].
\end{equation}
Now that the last exponential term is decoupled with respect to the replica indices, we obtain 
\begin{align}
  &\prod_{a, b, s \leq t} \delta\qty(NQ_{u,ab}^{st} - (\bm{u}_a^s)^\ten \bm{u}_b^t)\\
  =& \int \prod_{s \leq t}^T \dd \hat{q}_u^{st} \dd \hat{\chi}_u^{st} \exp \qty{ N n \sum_{t = 1}^T \beta_u^t \qty( \frac{q_u^{tt} \hat{q}_u^{tt} - \chi_u^{tt}\hat{\chi}_u^{tt}}{2} + \sum_{ s < t} (q_u^{st} \hat{q}_u^{st} - \chi_u^{st}\hat{\chi}_u^{st} ) + O(n)) } \nonumber \\
  &\times \prod_{i = 1}^N \qty[ \EE_{\bm{x}_u \sim \mathcal{N}(\bm{0}, \vb{\chi}_u)} \prod_{a = 1}^n \exp \sum_{t = 1}^T \beta_u^t \qty{ -\frac{\hat{q}_u^{tt}}{2} (u_{a,i}^t)^2 -  \sum_{s < t} \hat{q}_u^{st} u_{a,i}^s u_{a,i}^t + x_u^t u_{a,i}^t }].
\end{align}
Similiar expressions are derived for the state density term constraining $m_u^t$ and $R^t$ in \eqref{eq:state_density}, which is given by 
\begin{align}
  &\prod_{t =1}^T \prod_{a = 1}^n \delta\qty( N m_u^t - \sum_{i = 1}^N (u_{a,i}^t)^\ten u_{a,i}^t) \delta\qty( N R^t - \sum_{i = 1}^N (u_{a,i}^0)^\ten u_{a,i}^t) \\
  \propto &\int \prod_{t = 1}^T \dd \hat{m}_u^t \dd \hat{R}_u^t \exp \qty{ \sum_{t = 1}^T \beta_u^t \qty( -N n( m_u^t \hat{m}_u^t + R^t \hat{R}^t ) + \sum_{i = 1}^N  u_{a,i}^t  \qty(\hat{m}_u^t  u_i^\star + \hat{R}^t u_i^0)) }.
\end{align}
Thus, 
\begin{align}
  &\int \dd \bm{u}^0 \dd \bm{u}^\star  P(\bm{u}^0 | \bm{u}^\star) P(\bm{u}^\star) \prod_{a, t =1 }^{n, T} \dd \bm{u}^{t}_a  e^{ -\frac{\lambda \beta_u^t}{2} \norm{\bm{u}_a^t}_2^2   } \nonumber \\
  &\times \prod_{a \neq b, s \leq t} \delta\qty( N Q_{u, ab}^{st} - (\bm{u}^s_a)^\ten \bm{u}^t_b ) \prod_{a,t = 1}^{n,T} \delta\qty( N R^t_a - (\bm{u}^t_a)^\ten \bm{u}^0 ) \delta\qty( N m_{u, a}^t - (\bm{u}^t_a)^\ten \bm{u}^\star ) \\
  = &\int \prod_{s \leq t}^T \dd \hat{q}_u^{st} \dd \hat{\chi}_u^{st} \prod_{t = 1}^N \dd \hat{m}_u^t \dd \hat{R}^t \nonumber \\
  & \times \exp \qty{ N n \sum_{t = 1}^T \beta_u^t \qty( \frac{q_u^{tt} \hat{q}_u^{tt} - \chi_u^{tt}\hat{\chi}_u^{tt}}{2} - m_u^t \hat{m}_u^t - R^t \hat{R}^t + \sum_{ s < t} (q_u^{st} \hat{q}_u^{st} - \chi_u^{st}\hat{\chi}_u^{st} )) } \\
  &\times \qty[  \int \dd u^0  P(u^0 | u^\star) \dd u^\star P(u^\star) \EE_{\bm{x}_u} \qty{ \prod_{t = 1}^T \int \dd u^t e^{ \beta_u^t \qty[ -\frac{\hat{q}_u^{tt} + \lambda}{2} (u^t)^2 + \qty(x_u^t -\sum_{s < t} \hat{q}_u^{st} u^s + \hat{m}_u^t u^\star + \hat{R}^t u^0 )u^t ] } }^n]^N \nonumber.
\end{align}
Redefining $\hat{q}_u^{st}$ as $-\hat{q}_u^{st}$ for $s \neq t$, and using the saddle point approximation for large $N$ as well as a Laplace approximation for large $\beta_u^t$ in the the last equation above, 
one obtains the state density term for the $u$-variables as 
\begin{align}
  \begin{split}
 \exp nN \Extr \Bigg\{ &\sum_{t = 1}^T \beta_u^t  \qty[\frac{q_u^{tt} \hat{q}_u^{tt} - \chi_u^{tt} \hat{\chi}_u^{tt} }{2} - m_u^t \hat{m}_u^t - R^t \hat{R}^t - \sum_{s < t} ( q_u^{st} \hat{q}_u^{st} + \chi_u^{st} \hat{\chi}_u^{st} )] \\
 &- \EE_{\bm{x}_u} \min_{\bm{u}} \sum_{t = 1}^T \beta_u^t \qty[ \frac{\hat{q}_u^{tt} + \lambda}{2} (u^t)^2 - \qty(x_u^t + \sum_{s < t} \hat{q}_u^{st} u^s + \hat{m}_u^t u^\star + \hat{R}^t u^0 )u^t ]  \Bigg\}. 
  \end{split}
\end{align}
Keeping in mind that the limits of the $\beta_u$s are taken successively, one can notice that all terms $\qty{u^s}_{s < t}$ which appear under the optimization function with inverse temperature $\beta_u^t$ 
are all determined by the previous optimization functions with inverse temperatures $\beta_u^s$ for $s < t$. 
A sequence of random optimization problems can then be realized, resulting in the recursive structure described in the main text. 
One can also notice that the solution to the optimization problem is given by a Gaussian process \eqref{eq:gaussian_process}, 
resulting in the formula 
\begin{align}
  \exp nN \Extr \sum_{t = 1}^T \beta_u^t \Bigg\{ \frac{q_u^{tt} \hat{q}_u^{tt} - \chi_u^{tt} \hat{\chi}_u^{tt} }{2} - m_u^t \hat{m}_u^t - R^t \hat{R}^t - \sum_{s < t} ( q_u^{st} \hat{q}_u^{st} + \chi_u^{st} \hat{\chi}_u^{st} ) + \frac{\hat{q}_u^{tt} + \lambda}{2}\EE_{\bm{x}_u} \qty[ (\usf^t )^2] \Bigg\}.
 \end{align}
A similar computation follows for the $v$-variables, which offers 
\begin{align}
  \exp nN \Extr \sum_{t = 1}^T \beta_v^t \Bigg\{ \frac{q_v^{tt} \hat{q}_v^{tt} - \chi_v^{tt} \hat{\chi}_v^{tt} }{2} - m_v^t \hat{m}_v^t - \sum_{s < t} ( q_v^{st} \hat{q}_v^{st} + \chi_v^{st} \hat{\chi}_v^{st} ) + \frac{\hat{q}_v^{tt} + \lambda}{2}\EE_{\bm{x}_v} \qty[ (\vsf^t )^2] \Bigg\}.
\end{align}
 \subsection{Evaluation of the energy term}
 Here the object of interest is the energy term : 
 \begin{equation}
  \EE \Bigg[ \prod_{a = 1}^n  e^{-\beta_v^1 \ell ( h^0, k^1_a;y ) - \beta_u^1 l_{h^\star k^\star} ( h^1_a, k^1_a ) } \prod_{t =2 }^{T} e^{  - \beta_v^t \ell ( h^t_a, k^{t-1}_a;y ) - \beta_u^t \ell ( h^t_a, k^t_a;y ) } \Bigg],
 \end{equation}
 where the average $\EE$ is over the random fields $(h^\star, k^\star, h^0, \qty{h_a^t, k_a^t}_{a,t} )$, which are Gaussians distributed with a replica symmetric covariance: 
\begin{align}
  \begin{split}
  &\EE \qty[ (h^\star)^2 ] = \EE \qty[ (k^\star)^2 ] = \EE \qty[ (h^0)^2 ] = 1, \quad \EE \qty[ h^\star h^0 ] = m^0,\\
  &\EE \qty[ h^\star h_a^t ] = m_u^t, \quad \EE \qty[ h^0 h_a^t ] = R^t, \quad \EE \qty[ h_a^t h_b^s ] = q_u^{st} - (1-\delta_{ab}) \frac{\chi_u^{st}}{\beta^s_u}  \quad ( s \leq t, 1 \leq a,b \leq n ),\\
  &\EE \qty[ k^\star k_a^t ] = m_v^t,  \quad \EE \qty[ k_a^t k_b^s ] = q_v^{st} - (1-\delta_{ab}) \frac{\chi_v^{st}}{\beta^s_v},  \quad ( s \leq t, 1 \leq a,b \leq n ).
  \end{split}
\end{align}
Using the differential operator representation of the Gaussian average, i.e. 
\begin{equation}
  \EE_{\bm{x} \sim \mathcal{N}(\bm{0}, \vb{M})} \mathcal{F}(\bm{x}) =\eval{ \exp \qty(  \frac{1}{2} \sum_{i,j} M_{ij} \pdv{}{x_i}{x_j} ) \mathcal{F}(\bm{x})}_{\bm{x} = \bm{0}},
\end{equation}
the corresponding operator for this average is given by 
\begin{align}
\exp\Bigg[ \frac{1}{2} \dif^2_{h^\star} + m_0 \dif_{h^\star}\dif^2_{h^0} &+ \frac{1}{2} \dif_{k^\star}^2 + \sum_{t = 1}^{T} \qty( m_u^t \dif_{h^\star} + R^t \dif_{h^0} ) \sum_{a = 1}^n \dif_{h_a^t} \nonumber \\
&+  \frac{1}{2} \sum_{s , t = 1}^T q_u^{st} \qty( \sum_{a = 1}^n \dif_{h_a^t} )\qty( \sum_{a = 1}^n \dif_{h_a^s} ) + \frac{1}{2} \sum_{s,t = 1}^T \sum_{a = 1}^n \frac{\chi_u^{st}}{\beta_u^{\min(s,t)}} \dif_{h_a^t} \dif_{h_a^s}   \Bigg],
\end{align}
where we abuse notations as $\chi_u^{st} = \chi_u^{ts}$ and $ q_u^{st} = q_u^{ts}$. 
%Using again the multi-dimensional Hubbard Stratonovich transformation, we obtain
%\begin{align}
%  e^{ \frac{1}{2} \dif^2_{h^\star} + m_0 \dif_{h^\star} + \frac{1}{2}\dif^2_{h^0} } \EE_{\tilde{\bm{h}} \sim \mathcal{N}(\bm{0} , \vb{q}_u)} \prod_{a = 1}^n \exp \Bigg(\frac{1}{2} \sum_{s,t=1}^T \frac{\chi_u^{st}}{\beta_u^{\min(s,t)}} \dif_{h_a^t} \dif_{h_a^s} + \sum_{t = 1}^{T} \qty( \tilde{h}^t + m_u^t \dif_{h^\star} + R^t \dif_{h^0} )\dif_{h_a^t}    \Bigg) .
%  \end{align}
  We first simplify the last second-order differential operator in the above expression, which introduce Gaussian random variables with covariances $\qty{\chi_u^{st} / \beta_u^{\min(s,t)}}_{s,t}$. 
  Consider the Cholesky decomposition of this matrix, i.e. 
  \begin{equation}
    \mqty(  \chi_u^{11} / \beta_u^{1} & \chi_u^{12} / \beta_u^{1} &  \chi_u^{13} / \beta_u^1 & \cdots \\
    \chi_u^{21} / \beta_u^{1} & \chi_u^{22} / \beta_u^{2} &  \chi_u^{23} / \beta_u^2 & \cdots \\
    \chi_u^{31} / \beta_u^{1} & \chi_u^{32} / \beta_u^{2} &  \chi_u^{33} / \beta_u^3 & \cdots \\
    \vdots & \vdots & \vdots & \ddots ) = \mqty(L_{11} & 0 & 0 & \cdots \\ 
    L_{12} & L_{22} & 0 & \cdots \\ 
    L_{13} & L_{23} & L_{33} & \cdots \\
    \vdots & \vdots & \vdots & \ddots ) \mqty(L_{11} & L_{12} & L_{13} & \cdots \\
    0 & L_{22} & L_{23} & \cdots \\
    0 & 0 & L_{33} & \cdots \\
    \vdots & \vdots & \vdots & \ddots ),
  \end{equation}
  where $\vb{L}$ is a lower triangular matrix. It is not difficult to see that, in the limit successive $\beta_u^t$, 
  the decomposition is given by 
  \begin{equation}
    L_{st} = \frac{\chi_u^{st}}{\sqrt{ \beta_u^s \chi_u^{ss} }} + o( (\beta_u^s)^{-1/2} ), \qquad s \leq t. 
  \end{equation}
  Thus, the second-order differential operator can be expressed as
  \begin{gather}
    \exp\qty(\frac{1}{2} \sum_{s,t = 1}^T \frac{\chi_u^{st}}{\beta_u^{\min(s,t)}} \dif_{h_a^t} \dif_{h_a^s} ) = \exp\qty[ \frac{1}{2} \sum_{t=1}^T \qty( \sum_{s = t}^T L_{ts} \dif_{h_a^s} )^2 ]= \int \prod_{t =1}^T D z^t \exp\qty(- z^t  \sum_{s = t}^T L_{ts} \dif_{h_a^s} ) \nonumber \\
    =\int \prod_{t = 1}^T \dd z^t \exp \qty(  -\frac{\beta_u^t (z^t)^2}{\chi_u^{tt}} +   \sum_{s = 1}^t \frac{\chi_u^{st}}{\chi_u^{ss}} z^s \dif_{h_a^t} + o(\beta_u^t) ). 
  \end{gather}
  Note that the first-order differential operator is merely a classical translation operator. 
  The remaining second-order differential operators are handled by noticing the following identity: 
  \begin{gather}
    \sum_{i = 1}^d \dif_{x_i} \eval{ f(x_1, \cdots , x_d) }_{x_1 = \cdots = x_d = x} = \dif_x f(x, \cdots, x).
  \end{gather}
  This identity offers for a trial function $\mathcal{F},$
  \begin{align}
    &\exp \qty[ \sum_{t = 1}^{T} \qty( m_u^t \dif_{h^\star} + R^t \dif_{h^0} ) \sum_{a = 1}^n \dif_{h_a^t} +  \frac{1}{2} \sum_{s , t = 1}^T q_u^{st} \qty( \sum_{a = 1}^n \dif_{h_a^t} )\qty( \sum_{a = 1}^n \dif_{h_a^s}) ] \eval{\mathcal{F}\qty(\qty{h_a^t}_{a,t})}_{\qty{h_1^t = \cdots = h_n^t}_t } \nonumber \\
    =& \exp \qty[ \sum_{t = 1}^{T} \qty( m_u^t \dif_{h^\star} + R^t \dif_{h^0} ) \dif_{h^t} +  \frac{1}{2} \sum_{s , t = 1}^T q_u^{st} \dif_{h^t}\dif_{h^s} ] \mathcal{F}\qty(\qty{h^t}_{t}).
  \end{align}
  By applying all the same treatments to the $v$ (or $k$)-variables, the energy term is given by 
  \begin{align}
    \begin{split}
    &\exp \Bigg(  \frac{1}{2} \dif^2_{h^\star} + m_0 \dif_{h^\star} + \frac{1}{2}\dif^2_{h^0} + \frac{1}{2} \dif^2_{k^\star} \\
    + &  \sum_{t = 1}^T \qty( m_u^t \dif_{h^\star} \dif_{h^t}  + R^t \dif_{h^0} \dif_{h^t} + m_v^t \dif_{k^\star}\dif_{k^t}) + \frac{1}{2} \sum_{s,t = 1}^T \qty( q_u^{st} \dif_{h^t}\dif_{h^s} + q_v^{st} \dif_{k^t}\dif_{k^s} ) \Bigg)\\
    \times &\Bigg\{ \int \prod_{t = 1}^T \dd z^t \dd w^t \exp \qty(  -\frac{\beta_u^t(z^t)^2}{2\chi_u^{tt}} - \frac{\beta_v^t(w^t)^2}{2\chi_v^{tt}} +  \sum_{s = 1}^t \frac{\chi_u^{st}}{\chi_u^{ss}}  z^s \dif_{h^t} +  \sum_{s = 1}^t \frac{\chi_v^{st}}{\chi_v^{ss}}  w^s \dif_{k^t}  )\\
    \times &\exp\qty[-\beta_v^1 \ell ( h^0, k^1;y ) - \beta_u^1 \ell ( h^1, k^1;y ) ] \prod_{t =2 }^{T} \exp\qty[- \beta_v^t \ell( h^{t-1}, k^{t}; y ) - \beta_u^t \ell ( h^t, k^t;y ) ] \eval{\Bigg\}^n}_{ \substack{ \bm{h} = \bm{0}, \\ \bm{k} = \bm{0} } },\hspace{-24pt}
    \end{split}
  \end{align}
  where $y := h^\star k^\star$. 
  Retranslating the second-order differential operators as Gaussian averages and the first-order differential operators as translational operators, i.e. $e^{a \partial_x } f(x) = f(x+a)$, we finally find 
  \begin{align}
    1&+ n\EE \log \int \dd \bm{z} \dd \bm{w} \exp \Bigg\{ - \sum_{t = 1}^T \beta_v^t \qty[ \frac{(w^t)^2}{2\chi_v^{tt}} + \ell \qty( h^{t-1} + \sum_{s = 1}^{t-1} \frac{\chi_u^{s,t-1}}{\chi_u^{ss}} z^s , k^t +  \sum_{s = 1}^t \frac{\chi_v^{st}}{\chi_v^{ss}} w^s ;y ) ] \nonumber \\
    & - \sum_{t = 1}^T \beta_u^t \qty[ \frac{(z^t)^2}{2\chi_u^{tt}} + \ell \qty( h^t + \sum_{s = 1}^t \frac{\chi_u^{st}}{\chi_u^{ss}} z^s , k^t  +  \sum_{s = 1}^t \frac{\chi_v^{st}}{\chi_v^{ss}}w^s ; y )  ]  \Bigg\}  + O(n^2).
  \end{align}
  Obviously the $O(n)$ term is of interest, to which we apply the Laplace approximation for large $\beta_u, \beta_v$ to obtain 
  \begin{align}
    - n\EE \min_{\bm{z}, \bm{w}} \Bigg\{ &\sum_{t = 1}^T \beta_v^t \qty[ \frac{(w^t)^2}{2\chi_v^{tt}} + \ell \qty(z^{t-1} + h^{t-1} + \sum_{s = 1}^{t-2} \frac{\chi_u^{s,t-1}}{\chi_u^{ss}}z^s , w^t + k^t +  \sum_{s = 1}^{t-1} \frac{\chi_v^{st}}{\chi_v^{ss}} w^s ;y ) ] \nonumber \\
    & + \sum_{t = 1}^T \beta_u^t \qty[ \frac{(z^t)^2}{2\chi_u^{tt}} + \ell \qty( z^t + h^t + \sum_{s = 1}^{t-1} \frac{\chi_u^{st}}{\chi_u^{ss}} z^s , w^t +k^t  +  \sum_{s = 1}^{t-1} \frac{\chi_v^{st}}{\chi_v^{ss}} w^s ;y )  ]  \Bigg\} .
  \end{align}
One can again notice the same conditional structure as in the state density term; all variables $\qty{z^s, w^s}$ which has already appeared in a minimization problem 
  holding a higher inverse temperature can be considered as fixed. 
  Therefore, conditioned on the random variables $\{h^t, k^t\}$, for a minimization problem with inverse temperature $\beta_u^t$, $\qty{z^s}_{s < t}, \qty{w^s}_{s \leq t}$ can be considered a fixed variable, and thus minimization is only performed on a single variable $z^t$. 
  Likewise, for a minimization problem with inverse temperature $\beta_v^t$, $\qty{z^s}_{s < t}, \qty{w^s}_{s < t}$ can be considered as fixed, and thus minimization is only performed on the
  single variable $w^t$.
  This gives rise to the recursive structure described in the main text, \eqref{eq:effecive_solution}. 
\end{appendix}

\section{Proof of Claim 1}
\label{appendix:claim1_proof}

Here we prove that under finite $\kappa$ and zero initial overlap $m_0 = 0$, one cannot have a non-zero overlap $m_u^t$ for finite $t \geq 1$ in the limit of large $N$.
The proof is based on mathematical induction for the statement $m_u^t = m_v^t = \hat{m}_u^t = \hat{m}_v^t = 0$. 
The saddle point equation for $m_v^{t = 1}$ is given by 
\begin{equation}
   m_v^{1} = \EE_{0, \star} \qty[ \vsf^\star \vsf^1 ] = \frac{\hat{m}_v^1}{\hat{q}_v^{tt} + \lambda}. 
\end{equation}
Recalling that the dependency of $L_v^1$ on $h^\star$ only appears via $y = h^\star k^\star$, the saddle point equation for $\hat{m}_v^1$ is given by
\begin{equation}
    \hat{m}_v^{1} = -\kappa \EE_{h,k} \qty[ \frac{\dd^2}{\dd k^1 \dd k^\star} L_u^1 ] = -\kappa \EE_{h,k} \qty[ h^\star \frac{\dd^2}{ \dd k^1 \dd y}\ell(h^0, \phi_v^1 + w^1;y) ] = 0, 
\end{equation}
where we used that for an arbitrary differentiable function $f(x)$,
\begin{equation}
    \int Dz_1 Dz_2\ z_1 f(z_1z_2) = \int Dz_1 Dz_2 \dd y \ z_1 \delta(y - z_1z_2) f(y) = \int \frac{\dd z_1 \dd y}{2\pi} \  {\rm sgn}(z_1)e^{-\frac{z_1^2}{2}- \frac{y^2}{2z_1^2}} f(y) = 0. 
\end{equation}
This verifies $m_v^{t = 1} = \hat{m}_v^{t = 1} = 0.$ 
For $m_u^{t = 1}$ and $\hat{m}_u^{t = 1}$, we have 
   $ m_u^1 = \frac{\hat{m}_u^1}{\hat{q}_u^{tt} + \lambda},$ and 
\begin{equation}
    \hat{m}_u^1 = -\kappa \EE_{h,k} \qty[ \frac{\dd^2}{ \dd h^1 \dd h^\star} L_u^1 ] = -\kappa \EE_{h,k} \qty[ k^\star \frac{\dd^2}{ \dd h^0 \dd y} \ell( z^1 + \phi_u^1, w^1 + \phi_v^1;y) ].
\end{equation}
Recall that $\phi_v^1$ is only a function of $k^\star$ only through $y$, and $k^1$ is independent of $k^\star$ since $m_v^1 = 0$, and thus $z^1$ is also a function of $k^\star$ through $y$ only.
Thus, the same argument as above yields $\hat{m}_u^{t = 1} = m_u^{t = 1} = 0$.

Now, suppose that $m_u^s = m_v^s = \hat{m}_u^s = \hat{m}_v^s = 0$ for $s = 1, \cdots, t$. Thus $\qty{z^s, w^s, \phi_u^s, \phi_v^s}_{s \leq t}$ is only a function of $h^\star$ or $k^\star$ only through $y$, 
and $h^s, k^s$ is independent of $h^\star, k^\star$ for $s \leq t$. Then, 
\begin{equation}
    m_v^{t+1} = \EE_{0, \star} \qty[ \vsf^\star \vsf^{t+1} ] = \frac{\hat{m}_v^{t+1} + \sum_{s= 1}^{t} \hat{q}_u^{st} {m}_u^s }{\hat{q}_v^{t+1,t+1} + \lambda} =  \frac{\hat{m}_v^{t+1} }{\hat{q}_v^{t+1,t+1} + \lambda}.
\end{equation}
The saddle point equation for $\hat{m}_v^{t+1}$ is given by
\begin{equation}
    \hat{m}_v^{t+1} = -\kappa \EE_{h,k} \qty[ \frac{\dd^2}{\dd k^{t+1} \dd k^\star} L_u^{t+1} ] = -\kappa \EE_{h,k} \qty[ h^\star \frac{\dd^2 }{\dd k^{t+1} \dd y} \ell( z^t  + \phi_u^{t}, w^{t+1} + \phi_v^{t+1};y) ] = 0, 
\end{equation}
since $\displaystyle \frac{\dd^2}{\dd k^{t+1} \dd y} \ell( z^t  + \phi_u^{t}, w^{t+1} + \phi_v^{t+1};y)$ is a function of $h^\star$ and $k^\star$ only through $y = h^\star k^\star$. 
This yields $m_v^{t+1} = \hat{m}_v^{t+1} = 0$. 
As was the case of $t = 1$, the exact same arguments hold for $m_u^{t+1}$ and $\hat{m}_u^{t+1}$, which completes the proof.

\section{Numerical evaluation of the saddle point equations}\label{appendix:saddle}

The saddle point equations \eqref{eq:u_saddlepoint} and \eqref{eq:v_saddlepoint} must be solved numerically via fixed-point iteration, which is a non-trivial task due to the random averages in the expressions. 
However, due to the Gaussian nature of the stochastic process \eqref{eq:gaussian_process} the non-hatted variables can be 
calculated analytically. 
For instance, 
\begin{align}
  \label{eq:m_u_explicit}
    m_u^t &= \EE_{0,\star} \qty[ \usf^t \usf^\star ] = \frac{ \hat{m}^t + m_0\hat{R}^t + \sum_{s< t} \hat{q}_u^{st} m_u^s }{\hat{q}_u^{tt} + \lambda}, \\
    \label{eq:R_explicit}
    R^t   &= \EE_{0,\star} \qty[ \usf^t \usf^0 ] = \frac{ \hat{m}^t m_0 + \hat{R}^t  + \sum_{s< t} \hat{q}_u^{st} R^s }{\hat{q}_u^{tt} + \lambda}, \\
    \label{eq:q_u_explicit}
    q_u^{t^\prime t} &= \frac{1}{(\hat{q}_u^{tt} + \lambda)(\hat{q}_u^{\tpr \tpr} + \lambda)} \Big[ \hat{\chi}_u^{\tpr t} + \hat{m}_u^t \hat{m}_u^{\tpr} + \hat{R}^t \hat{R}^\tpr + (\hat{R}^t \hat{m}_u^\tpr + \hat{R}^\tpr \hat{m}_u^t) m_0 \nonumber \\
    & \qquad \qquad\qquad\qquad+ \sum_{s < t} \hat{q}_u^{st} (\hat{m}_u^\tpr m^s_u + \hat{R}^\tpr R^s ) + \sum_{\spr < \tpr} \hat{q}_u^{\spr\tpr} (\hat{m}_u^t m^\spr_u + \hat{R}^t R^\spr ) \nonumber \\
    & \qquad \qquad\qquad\qquad+ \sum_{s < t} \sum_{\spr < \tpr} \hat{q}_u^{st}\hat{q}_u^{\spr \tpr} q_u^{s \spr} + \sum_{s < t} \hat{q}_u^{st} \EE_{0, \star} \qty[ x_u^{\tpr} \usf^s ]+ \sum_{\spr < \tpr} \hat{q}_u^{\spr\tpr} \EE_{0, \star} \qty[x_u^{t} \usf^{\spr}] \Big],
\end{align}
where $ \Gamma^{ts}_u :=\EE_{0,\star}\qty[ x^t_u \usf^s]$ is given by the recursion
\begin{equation}
    \Gamma^{ts}_u = \frac{\chi_u^{st} + \sum_{\spr < s} \hat{q}_u^{s \spr} \Gamma^{t\spr} }{ \hat{q}_u^{ss} + \lambda}.
\end{equation}
In addition, 
\begin{align}
  \label{eq:chi_u_explicit}
  \chi_u^{st}& = \frac{\hat{q}_u^{t t} + \lambda}{2}  \pdv{}{\hat{\chi}_u^{st}} \EE_{0,\star} \qty[ (\usf^t)^2 ] =  \frac{\hat{q}_u^{t t} + \lambda}{2}  \EE_{0, \star} \qty[ \pdv[2]{(\usf^t)^2}{x_u^s}{x_u^t} ] \nonumber \\
  & = \EE_{0, \star} \qty[ \pdv{\usf^t}{x_u^s} ] =\frac{1}{\hat{q}_u^{tt} + \lambda} \qty( \delta_{st} + \sum_{\tpr < t} \hat{q}_u^{\tpr t} \EE_{0, \star} \qty[ \pdv{\usf^\tpr}{x_u^s} ]  ) = \frac{1}{\hat{q}_u^{tt} + \lambda} \qty( \delta_{st} + \sum_{\tpr < t} \hat{q}_u^{\tpr t} \chi_u^{s\tpr} ),
\end{align}
where the second equation follows from the explicit form of $\usf^t$, given in \eqref{eq:gaussian_process}. Note that $\chi_u^{st} = 0$ for $s > t$, 
which finally yields the result \eqref{eq:chi_u}. 
Given $\hat{\vb{\Theta}}_u^{t-1}$, the above expressions can then be calculated with $O(t^3)$ operations. Analogous expressions for the corresponding $v$ order parameters can also be obtained; 
in fact, one directly acquires those formulas by replacing the $u$-variables with the $v$-variables, and equating $R^t$ and $m_0$ to zero. 

The average with respect to $(h,k)$ can be evaluated numerically via Monte Carlo integration.
More explicitly, we prepare $N_{\rm MC}$ samples of $\{h^0, h^\star, k^\star\}$, which 
are then used to calculate $N_{\rm MC}$ samples of $\{ w^1, z^1, w^2, z^2, \cdots \}$ in this order. 
In this process, additional $N_{\rm MC}$ samples of $\{h^t, k^t\}_{t=1}^T$ are also generated 
according to the variance given in \eqref{eq:covariances}, whose elements are calculated consecutively in the saddle point equations \eqref{eq:u_saddlepoint} and \eqref{eq:v_saddlepoint}.  
The expectation is then calculated by averaging over these $N_{\rm MC}$ random samples.
To calculate the total second differentials of $L_u^t$ and $L_v^t$ in an efficient manner, 
we introduce the following auxillary functions and variables :

\begin{align}
    g_u^t(a,b) &:= \argmin_{z} \qty{ \frac{z^2}{2\chi_u^{tt}} + \ell( z + a, b;y ) }, \\
    g_v^t(a,b) &:= \argmin_{w} \qty{ \frac{w^2}{2\chi_v^{tt}} + \ell( a, w + b;y ) }, \\
    \gsf_u^t &:= g_u^t\qty(  h^t + \sum_{s = 1}^{t-1} \frac{\chi_u^{st}}{\chi_u^{ss}}z^s, k^t + \sum_{s = 1}^t \frac{\chi_v^{st}}{\chi_v^{tt}} w^s) , \\
    \gsf_v^t &:= g_v^t \qty( h^{t-1} + \sum_{s = 1}^{t-1} \frac{\chi_u^{s,t-1}}{\chi_u^{ss}}z^s, k^t + \sum_{s = 1}^{t-1} \frac{\chi_v^{st}}{\chi_v^{tt}} w^s ),\\
    \lsf_u^t &:= \ell \qty(  h^t + \sum_{s = 1}^{t} \frac{\chi_u^{st}}{\chi_u^{ss}}z^s, k^t + \sum_{s = 1}^t \frac{\chi_v^{st}}{\chi_v^{tt}} w^s;y ),\\
    \lsf_v^t &:= \ell \qty(  h^{t-1} + \sum_{s = 1}^{t-1} \frac{\chi_u^{s,t-1}}{\chi_u^{ss}}z^s, k^t + \sum_{s = 1}^{t} \frac{\chi_v^{st}}{\chi_v^{tt}} w^s ;y ).
\end{align}
Note that $z^t = \gsf_u^t$ and $w^t = \gsf_v^t$ in \eqref{eq:effecive_solution}, but they are defined as the value of $g_u^t$ and $g_v^t$ given two arguments. 
Without confusion, we also define the partial derivative $\partial_i \gsf_u^t, \partial_i \gsf_v^t, \partial_i \lsf_u^t, \partial_i \lsf_v^t$ 
as the partial derivative of $g_u^t, g_v^t, \ell(h^t + \cdots, k^t + \cdots;y), \ell(h^{t-1} + \cdots, k^t + \cdots;y)$ with respect to its $i (=1,2)$-th argument, respectively. 
The same partial derivatives with respect to variable $y$ are also defined as $\partial_y \gsf_u^t, \partial_y \gsf_v^t, \partial_y \lsf_u^t, \partial_y \lsf_v^t$.
Consider that 
\begin{align}
    \frac{\dd^2}{\dd h^t \dd h^{t^\prime}} L_u^t &=  \frac{\dd }{\dd h^\tpr} \partial_1 \ell( z^t + \phi_u^t, w^t + \phi_v^t;y )\\
    & = (\partial_1^2 \lsf_u^t )\qty(  \delta_{\tpr t} + \sum_{s = 1}^t\frac{\chi_u^{st}}{\chi_u^{ss}} \pdv{ z^s}{h^\tpr}  ) + (\partial_1 \partial_2 \lsf_u^t ) \sum_{s = 1}^t \frac{\chi_u^{st}}{\chi_u^{ss}}\pdv{w^s}{h^\tpr}, 
\end{align}
from a simple application of the chain rule. Define $A_\tpr^t := \pdv{z^t}{h^\tpr} $ and $B_\tpr^t := \pdv{w^t}{h^\tpr}$. These can be calculated recusively due to the following equation: 
\begin{align}
    A_\tpr^t := \pdv{z^t}{h^\tpr} &= \pdv{}{h^\tpr} g_u^t\qty(  h^t + \sum_{s = 1}^{t-1} \frac{\chi_u^{st}}{\chi_u^{ss}}z^s, k^t + \sum_{s = 1}^t \frac{\chi_v^{st}}{\chi_v^{ss}} w^s) \\
    & =(\partial_1 \gsf_u^t )\qty(\delta_{t\tpr} +  \sum_{s = 1}^{t-1} \frac{\chi_u^{st}}{\chi_u^{ss}} A_\tpr^s )+ (\partial_2 \gsf_u^t) \sum_{s = 1}^t \frac{\chi_v^st}{\chi_v^{ss}}B_\tpr^s\\
    B_\tpr^t := \pdv{w^t}{h^\tpr} &= \pdv{}{h^\tpr} g_v^t\qty(  h^{t-1} + \sum_{s = 1}^{t-1} \frac{\chi_u^{s,t-1}}{\chi_u^{ss}}z^s, k^t + \sum_{s = 1}^{t-1} \frac{\chi_v^{st}}{\chi_v^{ss}} w^s   ) \\
    &= (\partial_1 \gsf_v^t ) \qty( \delta_{t-1, \tpr} +  \sum_{s = 1}^{t-1} \frac{\chi_u^{s,t-1}}{\chi_u^{ss}} A_\tpr^s ) + (\partial_2 \gsf_v^t) \sum_{s = 1}^{t-1} \frac{\chi_v^{st}}{\chi_v^{ss}} B_\tpr^s.
\end{align}
Therefore, the variables $A_\tpr^t$ and $B_\tpr^t$ can be obtained by a bookkeeping procedure. 
The same argument holds for calculating  $\pdv{}{h^t}{h^\star} L_u^t$ and $\pdv{}{h^t}{h^0} L_u^t$, in which case we introduce the bookkeeping variables 
\begin{equation}
    A_\star^t := \pdv{z^t}{h^\star}, \quad A_0^t := \pdv{z^t}{h^0}, \quad B_\star^t := \pdv{w^t}{h^\star}, \quad B_0^t := \pdv{w^t}{h^0}. 
\end{equation}
Then, similar calculations yield 
\begin{align}
    A_\star^t &= k^\star \partial_y \gsf_u^t + (\partial_1 \gsf_u^t) \sum_{s = 1}^{t-1} \frac{\chi_u^{st}}{\chi_u^{ss}} A_\star^s + (\partial_2 \gsf_u^t) \sum_{s = 1}^t \frac{\chi_v^{st}}{\chi_v^{ss}} B_\star^s, \\
    B_\star^t &= k^\star \partial_y \gsf_v^t + (\partial_1 \gsf_v^t) \sum_{s = 1}^{t-1} \frac{\chi_u^{s,t-1}}{\chi_u^{ss}} A_\star^s + (\partial_2 \gsf_v^t) \sum_{s = 1}^{t-1} \frac{\chi_v^{st}}{\chi_v^{ss}} B_\star^s, \\
    A_0^t &= (\partial_1 \gsf_u^t) \sum_{s = 1}^{t-1} \frac{\chi_u^{st}}{\chi_u^{ss}} A_0^s + (\partial_2 \gsf_u^t) \sum_{s = 1}^t \frac{\chi_v^{st}}{\chi_v^{ss}} B_0^s, \\
    B_0^t &= (\partial_1 \gsf_v^t) \sum_{s = 1}^{t-1} \frac{\chi_u^{s,t-1}}{\chi_u^{ss}} A_0^s + (\partial_2 \gsf_v^t) \sum_{s = 1}^{t-1} \frac{\chi_v^{st}}{\chi_v^{ss}} B_0^s.
\end{align}
Using the book-keeping variables, the second derivatives of $L_u^t$ required to calculate \eqref{eq:u_saddlepoint} are provided via: 
\begin{align}
  \frac{\dd^2}{\dd h^t \dd h^\tpr} L_u^t &= (\partial_1^2 \lsf_u^t )\qty(  \delta_{\tpr t} + \sum_{s = 1}^t\frac{\chi_u^{st}}{\chi_u^{ss}} A_\tpr^s ) + (\partial_1 \partial_2 \lsf_u^t ) \sum_{s = 1}^t \frac{\chi_v^{st}}{\chi_v^{ss}} B_\tpr^s, \\
  \frac{\dd^2}{\dd h^t \dd h^\star} L_u^t &= k^\star \partial_y \partial_1 \lsf_u^t + (\partial_1^2 \lsf_u^t) \sum_{s = 1}^t \frac{\chi_u^{st}}{\chi_u^{ss}} A_\star^s + (\partial_1 \partial_2 \lsf_u^t) \sum_{s = 1}^t \frac{\chi_v^{st}}{\chi_v^{ss}} B_\star^s, \\
  \frac{\dd^2}{\dd h^t \dd h^0} L_u^t &= (\partial_1^2 \lsf_u^t) \sum_{s = 1}^t \frac{\chi_u^{st}}{\chi_u^{ss}} A_0^s + (\partial_1 \partial_2 \lsf_u^t) \sum_{s = 1}^t \frac{\chi_v^{st}}{\chi_v^{ss}} B_0^s.
\end{align}

The same calculations can be repeated for the $v$-variables, with slight modifications. 
Introducing the bookkeeping variables 
\begin{equation}
    C_\tpr^t = \pdv{z^t}{k^\tpr}, \quad C_\star^t = \pdv{z^t}{k^\star}, \quad D_\tpr^t = \pdv{w^t}{k^\tpr}, \quad D_\star^t = \pdv{w^t}{k^\star}, 
\end{equation}
one obtains the recursive equations 
\begin{align}
    C_\tpr^t &= 
    (\partial_1 \gsf_u^t ) \sum_{s = 1}^{t-1} \frac{\chi_u^{st}}{\chi_u^{ss}} C_\tpr^s + (\partial_2 \gsf_v^t) \qty( \delta_{t\tpr} + \sum_{s = 1}^t \frac{\chi_v^{st}}{\chi_v^{ss}}D_\tpr^s), \\
    D_\tpr^t &= 
    (\partial_1 \gsf_v^t ) \sum_{s = 1}^{t-1} \frac{\chi_u^{s,t-1}}{\chi_u^{ss}} C_\tpr^s  + (\partial_2 \gsf_v^t) \qty( \delta_{t\tpr} + \sum_{s = 1}^{t-1} \frac{\chi_v^{st}}{\chi_v^{ss}} D_\tpr^s), \\
    C_\star^t &= h^\star \partial_y \gsf_u^t + (\partial_1 \gsf_u^t) \sum_{s = 1}^{t-1} \frac{\chi_u^{st}}{\chi_u^{ss}} C_\star^s + (\partial_2 \gsf_u^t) \sum_{s = 1}^t \frac{\chi_v^{st}}{\chi_v^{ss}} D_\star^s, \\
    D_\star^t &= h^\star \partial_y \gsf_v^t + (\partial_1 \gsf_v^t) \sum_{s = 1}^{t-1} \frac{\chi_u^{s,t-1}}{\chi_u^{ss}} C_\star^s + (\partial_2 \gsf_v^t) \sum_{s = 1}^{t-1} \frac{\chi_v^{st}}{\chi_v^{ss}} D_\star^s.
\end{align}
Utilizing these variables, we have the convenient expression for the partial derivatives of $L_v^t$ given by 
\begin{align}
  \frac{\dd^2}{\dd k^t \dd k^\tpr} L_v^t &= ( \partial_2^2 \lsf_v^t ) \qty( \delta_{\tpr t} + \sum_{s = 1}^t \frac{\chi_v^{st}}{\chi_v^{ss}} D_\tpr^s ) + (\partial_1 \partial_2 \lsf_v^t ) \sum_{s = 1}^{t-1} \frac{\chi_u^{s,{t-1}}}{\chi_u^{ss}} C_\tpr^s, \\
  \frac{\dd^2}{\dd k^t \dd k^\star} L_v^t &= h^\star \partial_y \partial_2 \lsf_v^t + (\partial_2^2 \lsf_v^t) \sum_{s = 1}^{t} \frac{\chi_v^{st}}{\chi_v^{ss}} D_\star^s + (\partial_1 \partial_2 \lsf_v^t) \sum_{s = 1}^{t-1} \frac{\chi_u^{s,{t-1}}}{\chi_u^{ss}} C_\star^s. 
  %\frac{\dd^2}{\dd k^t \dd k^0} L_v^t &= (\partial_2^2 \lsf_v^t) \sum_{s = 1}^{t} \frac{\chi_v^{st}}{\chi_v^{ss}} D_0^s + (\partial_1 \partial_2 \lsf_v^t) \sum_{s = 1}^{t-1} \frac{\chi_u^{s,{t-1}}}{\chi_u^{ss}} C_0^s.
\end{align}

Calculating a single bookkeeping variable requires $O(t)$ operations, and thus the average of the partial derivatives of $L_u^t$ and $L_v^t$ can be calculated with $O(N_{\rm MC}t^2)$ operations. 
For all experiments, we employed $N_{\rm MC} = 10^8$ Monte Carlo samples. 
However, it is important to note that even with this substantial number of samples, 
calculating the trajectory of the order parameters can be susceptible to numerical instabilities.
This instability arises from the recursive bookkeeping process, where one must calculate the average of products of random variables. 
For instance, the update of $A_\tpr^t$ consists a sum of $A_{\tpr}^s$ for $s < t$, multiplied by $\partial_1 \gsf_u^t$, which is also a random variable in itself. 
Therefore, $A_\tpr^t$ consists of a composite product of $t$ random variables, which exhibit heavy-tailed behaviors. 
Consequently, a Monte Carlo approximation of their averages is susceptible to outliers in the samples. 
While we do employ the Lugosi-Mendelson estimator \cite{Lugosi19} to estimate the average in a robust manner, 
this does not completely eliminate the potential for numerical instabilities. 

\bibliography{ref}

%%%%%%%%% END TODO: CONTENTS

%%%%%%%%%% TODO: BIBLIOGRAPHY
% Provide your bibliography here. You have two options:

%%% FIRST OPTION
% Write your entries here directly, following the example below, including:
% Author(s), Title, Journal Ref. with year in parentheses at the end, followed by the DOI number.

%\begin{thebibliography}{99}
%\bibitem{1931_Bethe_ZP_71} H. A. Bethe, {\it Zur Theorie der Metalle. i. Eigenwerte und Eigenfunktionen der linearen Atomkette}, Zeit. f{\"u}r Phys. {\bf 71}, 205 (1931), \doi{10.1007\%2FBF01341708}.
%\bibitem{arXiv:1108.2700} P. Ginsparg, {\it It was twenty years ago today... }, \url{http://arxiv.org/abs/1108.2700}.
%\end{thebibliography}

%%% SECOND OPTION
% Use your bibtex library, formatted by the SciPost style file.
%\bibliography{SciPost_Example_BiBTeX_File.bib}

%%%%%%%%%% END TODO: BIBLIOGRAPHY

\end{document}